\documentclass[a4paper,10pt]{article}
\usepackage[utf8]{inputenc}
\usepackage[english]{babel}
\usepackage{latexsym,cite,mathrsfs}
\usepackage{amsfonts}
\usepackage{amssymb, amsthm}
\usepackage{xcolor}
\usepackage{marginnote}
\usepackage{amsmath}
\usepackage{a4wide}
\usepackage[affil-it]{authblk}

\numberwithin{equation}{section}

\usepackage{color,enumitem,graphicx}
\usepackage[colorlinks=true,urlcolor=blue,
citecolor=red,linkcolor=blue,linktocpage,pdfpagelabels,
bookmarksnumbered,bookmarksopen]{hyperref}
\usepackage{comment}

\usepackage[left=2.9cm,right=2.9cm,top=2.8cm,bottom=2.8cm]{geometry}
\usepackage[hyperpageref]{backref}

\makeatletter
\providecommand\@dotsep{5}
\def\listtodoname{List of Todos}
\def\listoftodos{\@starttoc{tdo}\listtodoname}
\makeatother

\usepackage{verbatim} 
\usepackage{a4wide}
\usepackage[hyperpageref]{backref}

\let\OLDthebibliography\thebibliography
\renewcommand\thebibliography[1]{
	\OLDthebibliography{#1}
	\setlength{\parskip}{1pt}
	\setlength{\itemsep}{1pt plus 0.3ex}
}

\usepackage[hyperpageref]{backref}

\usepackage{accents}

\newtheorem{thm}{Theorem}[section]
\newtheorem{lem}[thm]{Lemma}
\newtheorem{prop}[thm]{Proposition}
\newtheorem{cor}[thm]{Corollary}
\newtheorem{definition}[thm]{Definition}
\theoremstyle{definition}
\newtheorem{rem}[thm]{Remark}

\title{Rayleigh quotients of the level set  manifolds related to the nonlinear PDE}

\author{Yavdat Il'yasov\\

Institute of Mathematics, Ufa Federal Research Centre, RAS, 
	Chernyshevsky str. 112, 450008 Ufa, Russia\\
	Universidade Federal de Goi\'as, Instituto de Matem\'atica\\ 74690-900, Goi\^ania - GO - Brazil}
\begin{document}
\maketitle
\begin{abstract}
The main topic of this note is a discussion of applicability conditions of the Nehari manifold method depending on the value of parameters of equations. As the main tool, we apply the nonlinear generalized Rayleigh quotient method. 


\end{abstract}

\section{Introduction}

Let $W$ be a real Banach space, $\Phi_\lambda:W \to \mathbb{R}$ be a twice Fr\'{e}chet-differentiable functional, and  $\lambda$ be a real parameter. The Nehari manifold method introduced by Z. Nehari \cite{Neh1, Neh} in 1960 is a   powerful tool  for the investigation of equations of the variational form 
\begin{equation}
\label{I}
D\Phi_\lambda(u)=0,~~~u \in W, 
\end{equation}
which consists  of finding a extremal point $\hat{u}$ of the functional $\Phi_\lambda(u)$  subject  to  the \textit{Nehari manifold }    
$$
\mathcal{N}_{\lambda}=\{u\in W\setminus 0:~ ~\Phi'_\lambda(u):=D\Phi_\lambda(u)(u)=0\}. 
$$
This may lead to a solution to the problem \eqref{I}. Indeed,  assume that $\Phi''_\lambda(\hat{u}):=D^2\Phi_\lambda(\hat{u})(\hat{u},\hat{u})\neq 0$. Then under some general assumption the Lagrange multiplier rule due to {Lusternik} \cite{Lusternik}  yields $D\Phi_\lambda(\hat{u})+\mu D\Phi'_\lambda(\hat{u})=0$, for $\mu \in \mathbb{R}$, and consequently, $D\Phi_\lambda(\hat{u})=0$ since  $\Phi_\lambda'(\hat{u})=0$ and $\Phi''_\lambda(\hat{u})\neq 0$.

We call $\{\Phi''_\lambda(u)\neq 0,~\forall u \in \mathcal{N}_{\lambda}\}$ the applicability condition of the Nehari manifold method. Observe, $\mathcal{N}_{\lambda}$ is a zero level set of $\Phi'_\lambda(u)$, and thus, the applicability condition of the Nehari manifold method means that zero is a regular value of $\Phi'_\lambda(u)$.  In general, the question of the regularity (or lack thereof) of level sets of functionals associated with equations arises in many cases in the study of nonlinear problems. Below we will deal with finding the so-called solutions with prescribed energy $E$,  i.e., with solutions from the level set
$\{u \in W: \Phi_\lambda(u)=E\}$ (see  \cite{ilDiazR,ilStab}). Another type of level sets arises in the investigation of equations using the Pohozaev manifold (see \cite{ilDiaz2, ilegorov,ilTakac, PoS0}), and it seems this does not exhaust all the examples. 

The study of regularity (or lack thereof) of level sets  of functionals is based on the now classical theory (see, e.g., \cite{ lang, zeidl}). This problem takes on a different shade if the functional (equation) depends on a parameter. For example, the applicability condition of the Nehari manifold method may depends on the value of parameter $\lambda$. In this context, the question arises: \textit{how to find the values of the parameter $\lambda$ which corresponds to the regular zero level set of $\Phi'_\lambda$. } 
This leads us to the problem of finding the so-called \textit{Nehari manifold extreme value}, namely, limit points of the set of $\lambda$ where the applicability condition of the Nehari manifold method is satisfied. 

The main topic of this note is a discussion of the question: \textit{How to find the Nehari manifold extreme values?}

In  \cite{ilyaReil}, the so-called nonlinear generalized Rayleigh (NG-Rayleigh) quotient method has been introduced which allows one to find the Nehari manifold extreme values.  
The method is based on the analysis of NG-Rayleigh quotients whose critical values correspond to the Nehari manifold extreme values. 
To specify the principal idea of the method, 
let us consider  equation (\ref{I}) in the following particular form 
$$
D h(u)-\lambda Dg(u)=0.
$$
Assume for simplicity that $g'(u)\neq 0$, $\forall u\in W\setminus 0$. Testing the equation by $u \in W$ and  solving  it with respect to $\lambda$ we obtain the following parameter independent functional: 
\begin{equation}\label{RQ}
\mathcal{R}(u)=\frac{Dh(u)(u)}{Dg(u)(u)}=\frac{h'(u)}{g'(u)},~~~u\in W \setminus 0,
\end{equation}
which we call the Rayleigh quotient. 
Note that  $u$ belongs to $\mathcal{N}_\lambda$ if and only if it lies on the level set $\mathcal{R}(u)=\lambda$. Consider the one variable function $t \mapsto \mathcal{R}(tu)$,  $u \in W $, which we call the fibering function following Pohozaev \cite{Poh1,Poh}. Then 
\begin{equation*}
\mathcal{R}'(tu) :=\frac{d}{d t} \mathcal{R}(tu)=
\frac{1}{g'(tu)}\Phi_{\lambda}''(tu), ~~\forall u\in \mathcal{N}_\lambda, ~t>0.
\end{equation*}
Consequently, for $t>0$, $\mathcal{R}'(tu) \neq 0$, $\mathcal{R}(tu)=\lambda$ if and only if $\Phi''_{{\lambda}}(tu)\neq 0$, $tu \in \mathcal{N}_{\lambda}$. Thus, the finding of the applicable values $\lambda$ of the Nehari manifold method is reduced to the determining of the regular value of $\mathcal{R}(tu)$, $t>0$.   This reasoning  leads us to the following  idea.  

In many cases, the geometry of the fibering function $\mathcal{R}(tu)$,  $u \in W $ allows one   to determine a set of extreme points $t_i(u) \in \mathbb{R}^+\setminus 0$, $i=1,\ldots$ of $\mathcal{R}(tu)$. Under some general conditions, we can expect that the so-called  \textit{nonlinear generalized Rayleigh quotients} $\lambda_i(u):=\mathcal{R}(t_i(u )), ~~u \in W,~~ i=1,\ldots $ has  properties similar to that  the classical Rayleigh quotient has in the linear theory \cite{lanc}. The main idea of the NG-Rayleigh quotient method can be stated as follows: \, \textit{The set of the Nehari manifold extreme values can be found through  the set of critical values of the NG-Rayleigh quotients $(\lambda_i(u))$}. 

There is currently no general theory for this idea yet. Nevertheless, a number of its applications to specific problems allow us to conclude about its usefulness in solving different type of problems (see, e.g, \cite{ ilBob, bobkov, MarcCarlIl, ilDiazR,  ilStab, ilyaReil}). Below in this paper, are some of such examples of the application.

\begin{rem}
	In some cases, the Nehari manifold extreme values can also be found  using the so-called \textit{spectral analysis with respect to the fibering procedure} introduced in \cite{ilIzv, ilconcan}. In this approach, the Nehari manifold extreme values are found by solving the system 
	$$
	\left\{
	\begin{aligned}
		& h'(tu)-\lambda g'(tu)=0,\\
		& h''(tu)-\lambda g''(tu)=0,
		\end{aligned}
	\right.
	$$
	with respect to unknowns $t,\lambda$ (see, e.g., \cite{ilCherf, ilDiaz1, faraci, ilst, ilIzv, ilconcan, ilegorov, IlCritEx}).  
	However, in our opinion, the using of the NG-Rayleigh quotient method has a clear geometric meaning, which simplifies calculations and, ultimately, allows solving more complex problems.
\end{rem}

\begin{rem}
	Below we will see that the application of the NG-Rayleigh quotient method also makes it possible to constructively find useful variational formulations associated with equations. Research in this direction was motivated by the works  Pohozaev in \cite{Poh1,Poh}, where a fibering method was introduced to constructively find constrained minimization problems.
	 
\end{rem}

The paper is organized as follows. In Section \ref{sec:2}, we present the concept of the nonlinear Rayleigh quotient in an abstract setting. Section \ref{sec:3} is devoted to  an  example of the application of the NG-Rayleigh quotient method to a problem with convex-concave nonlinearity. In  particular, in this section, we discuses finding of the limit point for the branch of ground states obtained in the framework of the Nehari manifold method. In Section \ref{sec:4}, we present an example of the recursive application of the NG-Rayleigh quotient method to a problem where the fibering function $\Phi_{\bar{\lambda}}(t u)$  has  more than two critical points. 
Section \ref{sec:5} deals with the energy level Rayleigh quotient  which allow us to show the existence and nonexistence of solution with prescribed energy  for a zero mass  (zero frequency) problem.

\section{ Nonlinear generalized Rayleigh quotient}\label{sec:2}

We recall first some definitions and facts from the theory of manifolds. 
Let $W$ be a Banach space. Consider a Fr\'echet differentiable map $f: W \to \mathbb{R}$. We denote the Fr\'echet derivative  by $Df$.
For  $c \in V$,  the set
$f^{-1}(c)=\{u\in W: ~f(u)=c\}$ is said to be a $c$-level set of $f$.  
We call $u \in W$ the regular point of $f$ if  $Df(u):W \to \mathbb{R}$ is surjective; it is a critical point of $f$ otherwise. A point $c \in \mathbb{R}$  is said to be a regular value of $f$ if
every point of the level set $f^{-1}(c)$ is a regular point, and a critical value otherwise.

Let us consider the family of maps of the following form  $f_\lambda=h-\lambda g$,  where $h,g \in C^1(W; \mathbb{R})$, $\lambda \in \mathbb{R}$. For a given $c \in \mathbb{R}$, we are interested in the regularity  of $f^{-1}_\lambda(c)$ depends on the value of parameter $\lambda$. This problem can be investigate using  the following parameter independent functional
$$
\mathcal{R}(c;u):=\frac{h(u)-c}{g(u)}, ~~u \in W~ s.t. ~g(u)\neq 0,
$$
which we call   the \textit{Rayleigh quotient of the $c$-level manifold}  (\textit{Rayleigh quotient} for short) of $f_\lambda$. 

Assume that $\mathcal{R}(c;u)=\lambda$ for $\lambda \in \mathbb{R}$. Then  
$$
D\mathcal{R}(c;u)=\frac{1}{g(u)}(Dh(u)-\mathcal{R}(c;u) Dg(u))=\frac{1}{g(u)}Df_\lambda(u).
$$
Thus,   the map $Df_\lambda(u)$ is surjective if and only if $D\mathcal{R}(c;u)$ is surjective. From this we have
\begin{lem}\label{lem:lem1}
	   Suppose that $\lambda$ is a regular value of $ \mathcal{R}(c;\cdot)$,  and  ${f}^{-1}_\lambda(c)\cap\{u\in W:~g(u)= 0\} =\emptyset$. Then the level set
		${f}^{-1}_\lambda(c)$  is a $C^1$-manifold in $W$. Moreover $T_{u}(f^{-1}_\lambda(c))={\rm Ker} D \mathcal{R}(c;u)$   for every $u \in {f}^{-1}_\lambda(c)$ .  
\end{lem}
\begin{proof}
By the Regular Value Theorem  \cite{ lang, zeidl}
 the set $(\mathcal{R}(c;\cdot))^{-1}(\lambda)$ is a $C^1$-manifold and there holds $T_{u}((\mathcal{R}(c;\cdot))^{-1}(\lambda))={\rm Ker} D \mathcal{R}(c;u)$, $\forall u \in  (\mathcal{R}(c;\cdot))^{-1}(\lambda)$. Since ${f}^{-1}_\lambda(c)=(\mathcal{R}(c;\cdot))^{-1}(\lambda)$, this implies the proof.
\end{proof} 
 Define  
$$
\displaystyle{f_\lambda'(u):=\frac{d}{dt} (h(tu)-\lambda g(tu))|_{t=1}}\equiv Dh(u)(u)-\lambda D g(u)(u),~u \in W.
$$
The zero-level set
	$$
	{\mathcal N}_\lambda:=(f_\lambda')^{-1}(0)=\{u \in W:~f_\lambda'(u)=0\},
$$
is called a Nehari manifold  associated with $f_\lambda(u)$. 
A local minimum or maximum point of the function $f_\lambda$ subject to ${\mathcal N}_\lambda$ is called the extremal on the Nehari manifold.

\begin{lem}\label{lemNL}
 Let $\hat{u} \in {\mathcal N}_\lambda$ be an extremal point  of $f_\lambda$ on the Nehari manifold. Assume that $f_\lambda'(u)$ is Fr\'echet differentiable in an open neighborhood of $\hat{u} $ and $Df_\lambda'(\hat{u})$ is continuous at $\hat{u}$. Suppose that $f_\lambda''(\hat{u}):=Df_\lambda'(\hat{u})(\hat{u}) \neq 0$. Then $Df_\lambda(\hat{u})=0$.
	\end{lem}
\begin{proof} Due to the assumption we may apply the Lagrange multiplier rule  \cite{Lusternik} (see also Proposition 43.19. in \cite{zeidl}). Then, $Df_\lambda(\hat{u})+\mu Df'_\lambda(\hat{u})=0$, for some $\mu \in \mathbb{R}$. 
Testing this equality by $\hat{u}$ we obtain $\mu Df'_\lambda(\hat{u})(\hat{u})=\mu f_\lambda''(\hat{u})=0$. Since $f_\lambda''(\hat{u})\neq 0$, $\mu=0$, and therefore, $Df_\lambda(\hat{u})=0$.
\end{proof}

\begin{definition}
We call $\lambda $ the applicable value of the Nehari manifold method (applicable value of the Nehari manifold for short) if   $f_\lambda''(u)\neq 0$ for any $u \in {\mathcal N}_\lambda$. 	
\end{definition}
  Consider the  Rayleigh quotient of the Nehari manifold 
$$
\mathcal{R}^n(u)=\frac{h'(u)}{g'(u)},~ u \in W: g'(u)\neq 0.
$$
Observe, ${\mathcal N}_\lambda \neq \emptyset$, for any $\lambda \in {\rm Im}\,  \mathcal{R}^n$. Consider the fibering function
 $\mathcal{R}^n(tu)$ defined for each $u \in W$ s.t. $g(u)\neq 0$. We call $u \in W$ the fibering regular point of $\mathcal{R}^n(u)$ if $(\mathcal{R}^n)'(u) \neq 0$. A point $\lambda \in \mathbb{R}$  is said to be a fibering regular value of $\mathcal{R}^n$ if
every point of the level set $(\mathcal{R}^n)^{-1}(\lambda)$ is fibering regular. Evidently, the fibering regularity implies ordinary regularity in the above sense. We write $Z(\mathcal{R}^n)$ for the set of all fibering regular values of  $\mathcal{R}^n$.

\begin{cor}\label{corNL}  Assume that $f_\lambda, f_\lambda'\in C^1(W)$. 
  Suppose that $\lambda$ is a fibering regular value of $ \mathcal{R}^n$. Then $\lambda $ is the applicable value of the Nehari manifold method. Furthermore, 
	 if ${\mathcal N}_\lambda\cap\{u\in W:~g'(u)= 0\} =\emptyset$, $\forall \lambda \in Z(\mathcal{R}^n)$, then the set of fibering regular values of $\mathcal{R}^n$ coincides with the set of  applicable values of the Nehari manifold method.
	\end{cor}
	{\it Proof} is straightforward.

Let us show how this can be used by a simple example. Suppose that  $\forall u \in W$ s.t. $g'(u)\neq 0$:
\begin{quote}
	$(S)$:\, {\it $\mathcal{R}^n(tu)$ has no  critical  points in $\mathbb{R}^+\setminus 0$ except of global minimum or maximum points of $\mathcal{R}^n(tu)$.}
	\end{quote}
Introduce
\begin{align}
\lambda^{n,*}_{min}&=\sup_{u \in W:g(u)\neq 0}\,\lambda^{n}_{min}(u)\equiv \sup_{u \in W:g(u)\neq 0}\,\inf_{t \in \mathbb{R}^+}\, \mathcal{R}^n(tu), \label{lammin}\\
\lambda^{n,*}_{max}&=\inf_{u \in W:g(u)\neq 0}\,\lambda^n_{max}(u)\equiv\inf_{u \in W:g(u)\neq 0}\,\sup_{t \in \mathbb{R}^+}\, \mathcal{R}^n(tu). \label{lammax}
\end{align}

\begin{lem}\label{thm1}
  	Suppose   {\rm(S)}   and $-\infty\leq \lambda^{n,*}_{min} < \lambda^{n,*}_{max}\leq +\infty $. Assume that  ${\mathcal N}_\lambda\cap\{u\in W:~g'(u)= 0\} =\emptyset$, $\forall \lambda \in (\lambda^{n,*}_{min}, \lambda^{n,*}_{max})$. Then  $(\lambda^{n,*}_{min}, \lambda^{n,*}_{max})$ is an interval of applicability  of the Nehari manifold method.

\end{lem}
\begin{proof}  The proof immediately follows from Corollary \ref{corNL}.
\end{proof}   

The functionals  are given by
\begin{eqnarray}\label{NGRneh}
	\lambda_{min}(u):=\inf_{t \in \mathbb{R}^+}\, \mathcal{R}^n(tu),~~\lambda_{max}(u):=\sup_{t \in \mathbb{R}^+}\, \mathcal{R}^n(tu),
\end{eqnarray}
for $u \in W$ s.t. $g(u)\neq 0$ are called   \textit{the nonlinear generalized Rayleigh quotients} (\textit{the NG-Rayleigh quotients}). Notice that $\lambda_{min}(u)$, $\lambda_{max}(u)$ are $0$-homogeneous functionals, that is $\lambda_{min}(tu)=\lambda_{min}(u)$, $\lambda_{max}(tu)=\lambda_{max}(u)$, $\forall t>0$.
Observe that $0$-homogeneity is a basic property of the classical Rayleigh quotient used in the linear theory \cite{lanc}. 

\begin{rem}
Knowing the applicable values $\lambda$ of the Nehari manifold ${\mathcal N}_\lambda$ is important not only for finding solutions of equations by means of the Nehari minimization problem:
$$
\min\{f_\lambda(u):~ u \in {\mathcal N}_\lambda\}.
$$
The application of other methods, such as the critical point theory  or Ekeland's variational
principle requires that the Nehari manifold be $C^1$ (see, e.g., \cite{Bartsch, Bartsch2, Clark, Ekeland, ilconcan, szulkin, willem}).
\end{rem}
\begin{rem}
One of the advantages of using the Nehari manifold method is that it provides obtaining ground states of problems. In addition, in the case of parametric equations, this allows one to construct branches of solutions modulo values of the functional $f_\lambda$ \cite{ilBob, bobkov, ilst, ilIzv}, which is useful, for instance, in the investigation of the stability of ground states \cite{ilStab} and construction of solution at limit points \cite{bobkov,ilKaye}. 
\end{rem}
In general, for nonlinear problems, the critical points of
$\lambda_{min}(u)$, $\lambda_{max}(u)$ on $W$ do not necessarily correspond to the critical points of
$\Phi_{\lambda}(u)$. Indeed, suppose there exists $\hat{u} \in W\setminus 0$  such that  $D\lambda_{min}(\hat{u})=0$ in $W$. Since $\lambda_{min}(u)$ is $0$-homogeneous, we may assume that $\lambda_{min}(\hat{u})=\mathcal{R}^n(\hat{u})$. Hence  $D\mathcal{R}^n(\hat{u})=0$ which implies  $D^2\Phi_{\hat{\lambda}}(\hat{u})=0$ with $\hat{\lambda}=\lambda_{min}(\hat{u})$. Thus if $\hat{u}$ is a critical point of
$\Phi_{\lambda}(u)$ as well, then it has to satisfy simultaneously to the two equations $D\Phi_{\hat{\lambda}}(\hat{u})=0$, $D^2\Phi_{\hat{\lambda}}(\hat{u})=0$.
This is possible for linear equations, but impossible, in general, for  nonlinear problems.

However, there is a class of nonlinearly generalized Rayleigh quotients useful in finding solutions to nonlinear problems as well (see below and  \cite{MarcCarlIl, ilStab}), and which properties are similar to that of the classical  Rayleigh quotient.  Indeed, consider the   \textit{energy level Rayleigh quotient}  corresponding to the equation $Df_{\lambda}(u)\equiv Dh(u)-\lambda Dg(u)=0$, that is
$$
\mathcal{R}^e(E;u):=\frac{h(u)-E}{g(u)}, ~~u \in W: ~g(u)\neq 0.
$$
Here a value $E$ is called  energy (or action).  Observe the equality $D\mathcal{R}^e(E;u)=0$ implies that $Df_{\lambda}(u)=0$ and $f_{\lambda}(u)=E$, where $\lambda=\mathcal{R}^e(E;u)$, that is any critical point of $\mathcal{R}^e(E;u)$ corresponds to solution of the equation $Df_{\lambda}(u)=0$ with prescribed energy $E$.

For $u \in W$  s.t. $g(u)\neq 0$, consider the fibering function $\mathcal{R}^e(E;tu)$, $t>0$. Suppose that for every $ u \in W$ s.t. $g(u)\neq 0$, {\rm(S) is satisfied } with $\mathcal{R}^e(E;tu)$ instead of $\mathcal{R}^n(tu)$. Then as above we are able to introduce the following nonlinear generalized Rayleigh quotients
\begin{eqnarray*}\label{NGRneh2}
	\lambda_{min}^e(E;u):=\inf_{t \in \mathbb{R}^+}\, \mathcal{R}^e(E;tu),~~\lambda_{max}^e(E;u):=\sup_{t \in \mathbb{R}^+}\, \mathcal{R}^e(E;tu).
\end{eqnarray*} 
Hence, we have 
\begin{cor} 
$\lambda_{min}^e(E;u)$, $\lambda_{max}^e(E;u)$ are $0$-homogeneous functionals. Any critical point $\hat{u}$ of $\lambda_{min}^e(E;u)$ or $\lambda_{max}^e(E;u)$ is satisfied to equation $Df_{\lambda}(\hat{u})=0$ with prescribed energy level $f_{\lambda}(\hat{u})=E$ and frequency $\lambda=\lambda_{min}^e(E;\hat{u})$, $\lambda=\lambda_{max}^e(E;\hat{u})$, respectively.
\end{cor}

\noindent
\section{Extreme values of the convex-concave problem}\label{sec:3}
The complexity of  $\Phi_{\lambda}$ may be ranked depending on the number of critical points of the fibering functions $t\mapsto \Phi_{\lambda}(t u)$  $u \in W$. 
The simplest case is when $\Phi_{\lambda}(t u)$ has at most one critical point  for any $u \in W$ and $\lambda \in \mathbb{R}$, since in this case, in general, $\mathcal{N}_{\lambda}$ is a $C^1$-manifold. However, when $\Phi_{\lambda}(t u)$ may have two or more critical points  for some $u \in W$ and $\lambda \in \mathbb{R}$, the problem becomes more complicated, because of the Nehari manifold $\mathcal{N}_{\lambda}$ may contains a point $u$ where $\Phi''_{\lambda}(u)=0$.  This difficulty can be overcome by finding the corresponding Nehari manifold extreme values. In this section, we show an application of the nonlinear generalized Rayleigh quotient method to the case  when $\Phi_{\lambda}(t u)$ has at most two critical points.

Consider
\begin{equation}\label{eq:1}
\left\{
\begin{aligned}
	-&\Delta_p u  = \lambda |u|^{q-2}u+|u|^{\gamma-2}u &&\mbox{in}\ \ \Omega, \\
&u=0                                   &&\mbox{on}\ \ \partial\Omega.
\end{aligned}
\right.
\end{equation}
Here $\Delta_p u={\rm div}(|\nabla u|^{p-2}\nabla u)$ is the $p$-Laplacian, $\lambda\in \mathbb{R}$, $\Omega\subset \mathbb{R}^N$ is a bounded smooth domain, $1<q<p<\gamma<p^*$,  $p^*=\frac{pN}{N-p}$ if $p<N$, $p^*=+\infty$ if $p\geq N$. In what follows, the norm on the Sobolev space $W^{1,p}_0:=W^{1,p}_0(\Omega)$ we denote by $\|\cdot\|_1$. By a weak solution of \eqref{eq:1} we mean a critical point $u \in W^{1,p}_0$ of the energy functional 
$$
\Phi_\lambda(u):=\frac{1}{p}\int |\nabla u|^{p}\,dx -\lambda\frac{1}{q} \int |u|^{q}\,dx -\frac{1}{\gamma}\int |u|^{\gamma}\,dx.
$$
We construct solutions via the following two   Nehari manifolds minimization problems
\begin{align}
	&\hat{\Phi}^+_\lambda:=\min\{\Phi_\lambda(u):~~ u \in \mathcal{N}^+_\lambda\},\tag{${M}^+_\lambda$}\label{MinConv1}\\
	&\hat{\Phi}^-_\lambda:=\min\{\Phi_\lambda(u):~~u \in \mathcal{N}^-_\lambda\}.\tag{${M}^-_\lambda$}\label{MinConv2}
\end{align}
Here
\begin{align*}
	&\mathcal{N}^+_\lambda:=\{u \in W^{1,p}_0\setminus 0:~\Phi'_\lambda(u)=0,~\Phi''_\lambda(u)\geq 0\},\\
	&\mathcal{N}^-_\lambda:=\{u \in W^{1,p}_0\setminus 0:~\Phi'_\lambda(u)=0,~\Phi''_\lambda(u)\leq 0\},
\end{align*}
$\mathcal{N}_\lambda=\mathcal{N}^+_\lambda\cup \mathcal{N}^-_\lambda$ are the Nehari manifolds. A sequence $(u_m^\pm) \subset \mathcal{N}^\pm_\lambda$ is said to be minimizing of $(M^\pm_\lambda)$ if $\Phi_\lambda(u_m^\pm) \to \hat{\Phi}^\pm_\lambda$ as $m\to +\infty$. 

Observe that for $u \in W^{1,p}_0\setminus 0$, the fibering function $\Phi'_\lambda(tu)$ may have at most  two critical points: $t^+_\lambda(u)$, $t^-_\lambda(u)$ such that 
$0<t^+_\lambda(u)\leq t^-_\lambda(u)<+\infty$, $\Phi''_\lambda(t^+_\lambda(u)u)\geq 0$, $\Phi''_\lambda(t^-_\lambda(u)u)\leq 0$, and thus $t^\pm_\lambda(u)u\in \mathcal{N}^\pm_\lambda$. Furthermore, $t^+_\lambda(u)=t^-_\lambda(u)$ if and only if $\Phi''_\lambda(t^\pm_\lambda(u)u)=0$. Evidently,
\begin{equation}
\tag{${M}_\lambda$}\label{MinConvG}
	\hat{\Phi}^+_\lambda:=\min\{\Phi_\lambda(u):~~ u \in \mathcal{N}_\lambda\},
\end{equation}
and thus any minimizer of \eqref{MinConv1} is a ground state of \eqref{eq:1}.
By Lemma \ref{lemNL}, we have
 \begin{cor}\label{corrNACH}
	Let  $\lambda>0$. If a minimizer $\bar{u}_\lambda^\pm$ of  $(M^\pm_\lambda)$ satisfies $\Phi_\lambda''(\bar{u}_\lambda^\pm)\neq 0$, then $\bar{u}_\lambda^\pm$ is a weak solution of 
	\eqref{eq:1}. 
\end{cor}

Furthermore, we have
\begin{lem}\label{propCC}
Let $\lambda>0$. Any minimizing sequence of  $(M^\pm)$ has a nonzero limit point $u^\pm_\lambda \in W^{1,p}_0\setminus 0$ in the weak topology of $W^{1,p}_0$ and in the strong topology of $L^r$, $1<r<p^*$. Moreover, $\Phi_\lambda(u^\pm_\lambda)\leq \hat{\Phi}^\pm_\lambda$.
\end{lem}
 \begin{proof} The weak lower-semicontinuity of the norm of $W^{1,p}_0$ and the Sobolev theorem imply that $\Phi_\lambda(u)$ is weakly lower-semicontinuous on $W^{1,p}_0$. Furthermore, 
the functional $\Phi_\lambda(u)$ is coercive on $\mathcal{N}_\lambda$ since  by the Sobolev embedding we have
\begin{equation*}
	\Phi_\lambda(u)\geq  \frac{\gamma-p}{p\gamma}\|u\|^p_1 -\lambda\frac{\gamma-q}{q\gamma}\|u\|^q_{1},~~~\forall u \in \mathcal{N}_\lambda,~ \forall \lambda \in \mathbb{R},
\end{equation*}
where $p>q$.
Hence, for any $\lambda>0$, $\hat{\Phi}^\pm_\lambda>-\infty$, and any minimizing sequence of  $(M^\pm_\lambda)$   contains a subsequence $(u_m^\pm)$ which weakly in $W^{1,p}_0$ and strongly in $L^r$, $1<r<p^*$ converges to some limit point $u^\pm_\lambda$. This implies $\Phi'_\lambda(u^-_\lambda)=0$, $\Phi''_\lambda(u^-_\lambda)\leq 0$, and thus
$$
(p-q)\|u^-_\lambda\|^p_1\leq (\gamma-q)\|u^-_\lambda\|^\gamma_{L^\gamma}\leq C\|u^-_\lambda\|^\gamma_1,
$$
for some constant $C<+\infty$. Thus,  $u^-_\lambda\neq 0$. Note $\lim_{m\to+\infty}\Phi_\lambda(u_m^+)\geq 0$ if $u^+_\lambda = 0$. Hence, $u^+_\lambda\neq 0$ since $\Phi_\lambda(u_m^+)\to\hat{\Phi}^+_\lambda<0$. Evidently, $\Phi_\lambda(u^\pm_\lambda)\leq \hat{\Phi}^\pm_\lambda$.  
\end{proof}
To continue on, we need the following Rayleigh quotient 
\begin{equation*}
\mathcal{R}^n(u)=\frac{\int |\nabla u|^{p} dx-\int f |u|^{\gamma} dx}{\int |u|^{q} dx},~~u \in W^{1,p}_0\setminus 0.
\end{equation*}
Observe $\mathcal{N}_\lambda=\{u \in W^{1,p}_0: \mathcal{R}^n(u)=\lambda\}$.
The only critical point of $\mathcal{R}^n(su)$
is a   global maximum point of the function $\mathcal{R}^n(su)$ which can be found precisely:
\begin{equation*}
s_{max}(u)=\left(\frac{(p-q)\int|\nabla u|^{p} dx}{(\gamma-q)\int  |u|^{\gamma} dx}\right)^\frac{1}{\gamma-p}.  
\end{equation*}
Substituting  $s_{max}(u)$ into $\mathcal{R}^n(su)$ we obtain the following nonlinear generalized Rayleigh quotient: 
\begin{equation*}
\lambda(u):=
c_{p,q}\displaystyle{\frac{(\int|\nabla u|^{p} dx)^{\frac{\gamma-q}{\gamma-p}}}{\int |u|^{q} dx \,
	(\int |u|^{\gamma} dx)^{\frac{p-q}{\gamma-p}}}},  
\end{equation*}
where $c_{p,q}= \frac{\gamma-p}{p-q} \left( \frac{p-q}{\gamma-q}
\right)^{\frac{\gamma-q}{\gamma-p}} $. Notice that
\begin{equation}\label{lamEQUAL}
	\lambda'(u)=0,~\lambda(u)=\lambda~\Leftrightarrow~\Phi'_\lambda(u)=0,~\Phi''_\lambda(u)=0.
\end{equation}
Consider  the  nonlinear generalized Rayleigh  extremal value:    
\begin{equation}\label{LMaxN}
\lambda^*=\inf_{u \in W^{1,p}_0\setminus 0}\lambda(u)= 
c_{p,q}  \inf_{u \in W^{1,p}_0\setminus 0}\frac{(\int|\nabla u|^{p} dx)^{\frac{\gamma-q}{\gamma-p}}}{\int |u|^{q} dx \,
	(\int |u|^{\gamma} dx)^{\frac{p-q}{\gamma-p}} }. 
\end{equation} 
The Sobolev inequalities imply that $0<\lambda^*<+\infty$. 
Now we are able to prove the following criteria
\begin{lem}\label{propCCD}
Let $\lambda>0$. 
\begin{description}
		\item[$(1^o)$] If  $\lambda \leq \lambda(u^\pm_\lambda)$, then $u^\pm_\lambda$ is a nonzero minimizer of $(M^\pm_\lambda)$.
	\item[$(2^o)$] If $\lambda< \lambda(u^\pm_\lambda)$, then $u^\pm_\lambda$ is a nonzero weak solution of 
	\eqref{eq:1}.
\end{description}
\end{lem}
\begin{proof}
Let us show  $(1^o)$. The assumption $\lambda \leq \lambda(u^\pm_\lambda)$ implies the existence of $t^\pm_\lambda(u_\lambda^\pm)>0$. Hence, 
	\begin{align}\label{vsineq}
						\lambda=\mathcal{R}^n(t^\pm_\lambda(u_\lambda^\pm)u_\lambda^\pm)\leq\liminf_{m\to +\infty}\mathcal{R}^n(t^\pm(u_\lambda^\pm) u_m^\pm).
		\end{align}
	Since $\mathcal{R}^n( u_m^+)=\lambda$, $(\mathcal{R}^n)'( u_m^+)\geq 0$, $m=1,\ldots,$ \eqref{vsineq} yields that $1\leq t^+(u_\lambda^+)$, and thus 
	$$
	\Phi_\lambda(t^+(u_\lambda^+)u_\lambda^+)\leq \Phi_\lambda(u_\lambda^+)\leq \liminf_{m\to +\infty}\Phi_\lambda(u_m^+)=\hat{\Phi}^+_\lambda. 
	$$
	In view of that $t^+(u_\lambda^+)u_\lambda^+ \in \mathcal{N}^+_\lambda$, this implies that  $u_\lambda^+$ is a minimizer of \eqref{MinConv1}.
Similarly, \eqref{vsineq} implies that $ \Phi_\lambda(t^-(u_\lambda^-) u_m^-)\leq \Phi_\lambda(u_m^-)$, $m=1,\ldots$, and thus 
	$$
	\Phi_\lambda(t^-(u_\lambda^-)u_\lambda^-) \leq\liminf_{m\to +\infty}\Phi_\lambda(t^-(u_\lambda^-) u_m^-)\leq \liminf_{m\to +\infty}\Phi_\lambda(u_m^-)=\hat{\Phi}^-_\lambda, 
	$$
	which yields that  $u_\lambda^-$ is a minimizer of  \eqref{MinConv2}. From here and Corollary \ref{corrNACH} it follows
	$(2^o)$.	
\end{proof}

\begin{thm}\label{cor444}
Assume that  $\lambda\in (0,\lambda^*)$. Then \eqref{eq:1} has two distinct weak positive  solutions $u_\lambda^+$, $u_\lambda^-$ such that $\Phi_\lambda(u_\lambda^+)> 0$, $\Phi_\lambda''(u_\lambda^-)< 0$. Moreover,  $u_\lambda^+$ is a ground state of \eqref{eq:1}, $u_\lambda^\pm \in C^{1,\alpha}(\overline{\Omega})$, $\alpha \in (0,1)$.	
\end{thm}
\begin{proof}
 Observe for $\lambda\in (0,\lambda^*)$, we have $\lambda< \lambda^*\leq \lambda(u^\pm_\lambda)$. Hence Lemma  \ref{propCCD}  implies that  $u_\lambda^\pm$  is a weak solution of \eqref{eq:1}.  Since $\Phi_\lambda(|u|)=\Phi_\lambda(u)$ for $u\in W^{1,p}_0$, we may assume that $ u^\pm_\lambda\geq 0$ in $\Omega$. 
By the bootstrap argument	and the Sobolev embedding theorem
it follows that $u^\pm_\lambda \in L^\infty$. Then $C^{1,\alpha}$ -regularity results of DiBenedetto \cite{DiBened} and Tolksdorf \cite{Tolksdorf} (interior regularity) combined with  Lieberman \cite{Lieber} (regularity up to
the boundary) yield $u_\lambda^\pm \in C^{1,\alpha}(\overline{\Omega})$, for some $\alpha \in (0,1)$. Finally, since $p<\gamma$, the Harnack
inequality due to Trudinger \cite{trudin} implies that $u_\lambda^\pm>0$ in $\Omega$.
\end{proof}

We need the following lemma on a continuation to a limit point $\lambda^o>0$.
\begin{lem}\label{lemCC2}
	Assume that  $\lambda^o>0$. Let $\lambda_m \to \lambda^o$ as $m \to +\infty$, and $(u_{\lambda_m}^\pm)$ is a sequence of minimizers of $(M^\pm_{\lambda_m})$ such that  $D\Phi_{\lambda_m}(u_{\lambda_m}^\pm)=0$, $m=1,\ldots$.   Then there exists a weak solution $u_{\lambda^o}^\pm$ of \eqref{eq:1}, which is  a minimizer  of $(M^\pm_{\lambda^o})$. Moreover,  $\hat{\Phi}_{\lambda_m}^\pm\to \hat{\Phi}_{\lambda^o}^\pm$ as $m \to +\infty$. 
\end{lem}
\begin{proof}
	Let $\lambda_m \to \lambda^o$ as $m \to +\infty$, and $(u_{\lambda_m}^\pm)$ is a sequence of minimizers of $(M^\pm_{\lambda_m})$ such that  $D\Phi_{\lambda_m}(u_{\lambda_m}^\pm)=0$, $m=1,\ldots$.  Analysis similar to that in the proof of Lemma \ref{propCC} shows that there exist nonzero weak limit point $u^\pm_{\lambda^o} \in W^{1,p}_0\setminus 0$ of  $(u_{\lambda_m}^\pm)$. Hence passing to the limit in $D\Phi_{\lambda_m}(u_{\lambda_m}^\pm)=0$ we obtain $D\Phi_{\lambda^o}(u_{\lambda^o}^\pm)=0$. Thus, $u_{\lambda^o}^\pm$ is a weak solution of \eqref{eq:1}. 
	
Moreover, since $0=\Phi_{\lambda^o}'(u_{\lambda^o}^\pm)\leq \liminf_{m \to +\infty}\Phi_{\lambda_m}'(u_{\lambda_m}^\pm)=0$, it follows that $u_{\lambda_m}^\pm \to 
u_{\lambda^o}^\pm$ strongly in $W^{1,p}_0$, and, in particular, $\hat{\Phi}_{\lambda_m}^\pm\to \Phi_{\lambda^o}(u_{\lambda^o}^\pm)$ as $m\to +\infty$.

To prove that $u_{\lambda^o}^\pm$ is a minimize of $(M^\pm_{\lambda^o})$, it is sufficient to show that  $\hat{\Phi}_{\lambda^o}^\pm=\Phi_{\lambda^o}(u_{\lambda^o}^\pm)$. Conversely, suppose that $\hat{\Phi}_{\lambda^o}^\pm<\Phi_{\lambda^o}(u_{\lambda^o}^\pm)$. Then there exists $w^\pm \in \mathcal{N}_{\lambda^o}^\pm$ such that
$\Phi_{\lambda^o}(w^\pm)=\Phi_{\lambda^o}(u_{\lambda^o}^\pm)-\kappa$ for some $\kappa>0$. Evidently, for any $\epsilon>0$, one can find $m_\epsilon$ such that $|\Phi_{\lambda^o}(w^\pm)-\Phi_{\lambda_m}(t^\pm_{\lambda_m}(w^\pm)w^\pm)|<\epsilon$ and $|\hat{\Phi}_{\lambda_m}^\pm-\Phi_{\lambda^o}(u_{\lambda^o}^\pm)|<\epsilon$ for every $m>m_\epsilon$. Hence we have
$$
\Phi_{\lambda^o}(u_{\lambda^o}^\pm)-\epsilon<\hat{\Phi}_{\lambda_m}^\pm\leq\Phi_{\lambda_m}(t^\pm_{\lambda_m}(w^\pm)w^\pm) <\Phi_{\lambda^o}(w^\pm)+\epsilon=\Phi_{\lambda^o}(u_{\lambda^o}^\pm)-\kappa+\epsilon
$$
which implies a contradiction. Thus $\hat{\Phi}_{\lambda^o}^\pm=\Phi_{\lambda^o}(u_{\lambda^o}^\pm)$, and $\hat{\Phi}_{\lambda_m}^\pm\to \hat{\Phi}_{\lambda^o}^\pm$ as $m \to +\infty$.
\end{proof}
We need aslo the  zero-level energy Rayleigh quotient 
	\begin{equation*}
	\mathcal{R}^e(u):=\mathcal{R}^e(0;u)=\frac{\frac{1}{p}\int |\nabla u|^{p} dx-\frac{1}{\gamma}\int f |u|^{\gamma} dx}{\frac{1}{q}\int |u|^{q} dx},~~u \in W^{1,p}_0\setminus 0.
	\end{equation*}
	Notice that $\mathcal{R}^e(u)=\lambda~\Leftrightarrow ~ \Phi_\lambda(u)=0$. A computation similar to that has been used above for $\mathcal{R}^n(u)$ shows that  the unique critical point of $\mathcal{R}^e(su)$ in $s>0$ is a global maximum point $s^e_{max}(u)$ and one can introduce the corresponding nonlinear generalized Rayleigh quotient $\lambda^e(u)=\mathcal{R}^e(s^e_{max}(u)u)$. Moreover, for $u \in W^{1,p}_0\setminus 0$, 
\begin{align}
	&D\lambda^e(u)=0,~ \lambda^e(u)=\lambda^e~\Leftrightarrow ~D\Phi_{\lambda^e}(s^e_{max}(u)u)=0, ~\Phi_{\lambda^e}(s^e_{max}(u)u)=0\\
	& \lambda^e(u)=c_{pq}\lambda(u),~\mbox{where}~ c_{pq}=qp^{\frac{(p-q)}{(\gamma-p)}}/p^\frac{(\gamma-q)}{(\gamma-p)}. 
\end{align}
It is not hard to show that  ${\lambda}^e(u)<\lambda(u)$, $\forall u \in W^{1,p}_0\setminus 0$.
\begin{thm}\label{thmCC2}
	For $\lambda=\lambda^*$, \eqref{eq:1} has two distinct weak positive  solutions $u_{\lambda^*}^+$, $u_{\lambda^*}^-$ such that $\Phi_{\lambda^*}''(u_{\lambda^*}^+)> 0$, $\Phi_{\lambda^*}''(u_{\lambda^*}^-)< 0$, $\hat{\Phi}^+_{\lambda^*}<\hat{\Phi}^-_{\lambda^*}$. Moreover, $u_{\lambda^*}^+$, $u_{\lambda^*}^-$ are minimizers  of \eqref{MinConv1}, \eqref{MinConv2}, respectively; $u_{\lambda^*}^+$ is a ground state of \eqref{eq:1}; $u_\lambda^\pm \in C^{1,\alpha}(\overline{\Omega})$, $\alpha \in (0,1)$.
\end{thm}
\begin{proof} The existence of the weak  solutions $u_{\lambda^*}^+$, $u_{\lambda^*}^-$ which are minimizers  of \eqref{MinConv1}, \eqref{MinConv2}, respectively,
follows from Theorem \ref{cor444} and Lemma \ref{lemCC2}. 

Suppose, contrary to our claim, that $\Phi_{\lambda^*}''(u_{\lambda^*}^\pm)=0$. Then by \eqref{lamEQUAL}, $\lambda^*=\lambda(u_{\lambda^*}^\pm)$. Since \eqref{LMaxN}, this implies that $D\lambda(u_{\lambda^*}^\pm)=0$, and consequently $D\lambda^e(u_{\lambda^*}^\pm)=0$. Hence   $D\Phi_{{\lambda}^e}(s^e_{max}(u_{\lambda^*}^\pm)u_{\lambda^*}^\pm)=0$, where ${\lambda}^e:=\lambda^e(u_{\lambda^*}^\pm)$. Since $D\Phi_{{\lambda}^*}(u_{\lambda^*}^\pm)=0$ and ${\lambda}^e<\lambda^*$,  we get a contradiction.	The rest of the proof runs as in Theorem \ref{cor444}.
		
\end{proof}

\par
\noindent
In the case $\lambda>\lambda^*$, we have
	
\begin{thm}\label{thmCC3}
There exists $\bar{\lambda}\in (\lambda^*, +\infty]$ such that 	for any $\lambda\in (\lambda^*,\bar{\lambda})$ problem  \eqref{eq:1} has a ground state $u_\lambda^+$ such that $\Phi''_\lambda(u_\lambda^+)> 0$. Moreover $u_\lambda^+$ is a minimizer of   \ref{MinConvG}, $u_\lambda^+ \in C^{1,\alpha}(\overline{\Omega})$, $\alpha \in (0,1)$, $u_\lambda^+>0$ in $\Omega$. 
\end{thm}
\begin{proof} 
By Lemma \ref{propCCD}, it is sufficient to show that  there exists $\bar{\lambda}>\lambda^*$ such that $\lambda<\lambda(u_\lambda^+)$,	for any $\lambda\in (\lambda^*,\bar{\lambda})$.  

Suppose this is false. Then one can find a sequence $\lambda_m>\lambda^*$, $m=1,\ldots$ such that $\lambda_m \to \lambda^*$ as $m\to +\infty$ and $\lambda(u_{\lambda_m}^+)\leq \lambda_m$, $\forall m=1,\ldots$. Analysis similar to that in the proof of Lemma \ref{propCC} shows that there exists a nonzero weak limit point $\tilde{u}^+ \in W^{1,p}_0\setminus 0$ and a subsequence, which we still denote by
$(u_{\lambda_m}^+)$, such that $u_{\lambda_m}^+ \to \tilde{u}^+$ weakly in $W^{1,p}_0$ and strongly in $L^\beta$, $1<\beta<2^*$. Moreover, we may assume that  $\lim_{m\to +\infty}\hat{\Phi}^+_{\lambda_m}= \tilde{\phi}^+$  for some $\tilde{\phi}^+\leq 0$.

Suppose that the  convergence $u_{\lambda_m}^+ \to \tilde{u}^+$ is not strong in $W^{1,p}_0$. Since $\lambda^*\leq \lambda(\tilde{u}^+)$, 
  there exists $t_{\lambda^*}^+(\tilde{u}^+)>0$ such that 
$$
\lambda^*=\mathcal{R}^n(t_{\lambda^*}^+(\tilde{u}^+)\tilde{u}^+)<\liminf_{m\to +\infty}\mathcal{R}^n(t_{\lambda^*}^+(\tilde{u}^+)u_{\lambda_m}^+)\leq \liminf_{m\to +\infty}\lambda(u_{\lambda_m}^+)\leq \liminf_{m\to +\infty}\lambda_m= \lambda^*.
$$
We get a contradiction, and thus  $u_{\lambda_m}^+ \to \tilde{u}^+$ strongly in $W^{1,p}_0$. Consequently, $\Phi_{\lambda_m}(u_{\lambda_m}^+) \to \Phi_{\lambda^*}(\tilde{u}^+)= \tilde{\phi}^+$ as $m\to +\infty$, and  $\tilde{u}^+ \in \mathcal{N}^+_{\lambda^*}$
  
	Observe, $\Phi_{\lambda^*}(\tilde{u}^+)\leq \hat{\Phi}^+_{\lambda^*}$. Indeed, by Theorem \ref{thmCC2},
	there exists a minimizer ${u}_{\lambda^*}^+ $ of $(M^+_{\lambda^*})$ and $\lambda^*<\lambda(u_\lambda^*)$. Then for sufficiently large $m$, $\lambda_m<\lambda(u_\lambda^*)$, and thus there exists $t_{\lambda_m}^+(u_{\lambda^*}^+)>0$. Evidently, $t_{\lambda_m}^+(u_{\lambda^*}^+)\to 1$ as $m\to +\infty$. Hence,   
$$
\Phi_{\lambda^*}(\tilde{u}^+)=\tilde{\phi}^+=\lim_{m\to +\infty}\hat{\Phi}^+_{\lambda_m}\leq \lim_{m\to +\infty}\Phi_{\lambda_m}(t_{\lambda_m}^+(u_{\lambda^*}^+)u_{\lambda^*}^+)\leq  \lim_{m\to +\infty}\Phi_{\lambda^*}(t_{\lambda_m}^+(u_{\lambda^*}^+)u_{\lambda^*}^+)= \hat{\Phi}^+_{\lambda^*}.
$$
Thus,  $\tilde{u}^+$ is a minimizer of $(M^+_{\lambda^*})$.  Observe, $ \lambda(\tilde{u}^+)=\lim_{m\to +\infty}\lambda(u_{\lambda_m}^+)\leq\lambda^*$ since $\lambda(u_{\lambda_m}^+)\leq \lambda_m$,  $\forall m=1,\ldots$. Since $\mathcal{R}^n(\tilde{u}^+)=\lambda^*$, this implies $\lambda(\tilde{u}^+)=\lambda^*$, and consequently, $\Phi_{\lambda^*}''(\tilde{u}^+)=0$. Hence, we have $\tilde{u}^+ \in  \mathcal{N}^+_{\lambda^*}\cap  \mathcal{N}^-_{\lambda^*}$.  Using Theorem \ref{thmCC2}, we conclude that $\Phi_{\lambda^*}(\tilde{u}^+)=\hat{\Phi}^+_{\lambda^*}<\hat{\Phi}^-_{\lambda^*}\leq \Phi_{\lambda^*}(\tilde{u}^+)$, which is a contradiction. The rest of the proof runs as in Theorem \ref{cor444}.
\end{proof}

One can make the following conclusion about a limit point of the branch of ground states obtained by Nehari manifold method.

For $\lambda>0$, denote by $
	G_\lambda$
the set of ground states of \eqref{eq:1}, and by 
$$
	G^n_\lambda:=\{u \in W^{1,p}_0\setminus 0: \hat{\Phi}^+_\lambda= \Phi_\lambda(u),~\lambda<\lambda(u)\}
	$$ 
the subset of minimizes of \eqref{MinConvG}, and $\overline{G^n_\lambda}:=\{u \in W^{1,p}_0\setminus 0: \hat{\Phi}^+_\lambda= \Phi_\lambda(u),~\lambda\leq \lambda(u)\}$. The above results yield that for any $\lambda>0$, $G^n_\lambda=G_\lambda$  if 	$G^n_\lambda \neq \emptyset$. On the other hand, by in large, it is possible that $G_\lambda \neq \emptyset$ whereas
$G^n_\lambda = \emptyset$. 	Thus, it makes sense to consider the ground states obtained by means of Nehari manifold method as a particular branch. The family of set $G^n_\lambda$ which satisfies $G^n_\lambda \neq \emptyset$, we call  the   Nehari manifold ground states branch. 
	 
It is known  that the convex-concave problem \eqref{eq:1} possesses upper bound of the set positive solutions, namely, there exists $\lambda^b\in (0, +\infty)$ such that for any $\lambda>{\lambda^b}$, \eqref{eq:1} has no positive solutions. For the proof of this assertion, in the case $p=2$, we refer the reader to Ambrosetti, Brezis, Cerami \cite{AmBrCer}. The case $p\neq 2$, can be handled in much the same way using Picone's identity \cite{Aleg, IlCras}.

\begin{thm}\label{thmF}
There exists a limit value $\lambda^f\in (\lambda^*, +\infty)$ of the branch of the Nehari manifold ground states  such that:	(i) for any $\lambda\in (0,\lambda^f)$, $G^n_\lambda\neq \emptyset$;
(ii) $G_{\lambda^f}=\overline{G^n_{\lambda^f}}\neq \emptyset$; (iii)	there exists a sequence 
	$(\lambda_m)_{m=1}^\infty \subset (\lambda^f,+\infty)$ such that $\lambda_m \downarrow \lambda^f$,   and  $G^n_{\lambda_m}= \emptyset$, $\forall m=1,\ldots$.
 
\end{thm}
\begin{proof} The proof of the first part follows from Theorems \ref{cor444}, \ref{thmCC2}, \ref{thmCC3}. Suppose that there exists $\lambda\in (0,+\infty)$ such that $G^n_\lambda=\emptyset$. Then $\lambda^f=\sup\{\lambda>0:~G^n_\lambda\neq \emptyset\}<+\infty$. By the above results, $\lambda^f>\lambda^*$. Lemma \ref{lemCC2} yields that  there exists  a limit solution $u_{\lambda^f}^+$  of \eqref{eq:1} which is a minimizer  of $(M^+_{\lambda^f})$, and thus $u_{\lambda^f}^+ \in G_\lambda$. The proof of (iii) is straightforward.
\end{proof}

\begin{rem}
 We anticipate that $\lambda^f$ coincides with $\lambda^b$.
\end{rem}

\begin{rem} The existence of two distinct positive solutions of \eqref{eq:1} can be obtained  by the Nehari manifold method	 without finding the Nehari manifold extreme value $\lambda^*$, if $\lambda$ is sufficiently close to $0$,  where it can be shown by appropriate estimates that  $\mathcal{N}^0_\lambda:=\{u \in W^{1,p}_0\setminus 0:~\Phi'_\lambda(u)=0,~\Phi''_\lambda(u)= 0\}=\emptyset$ (see, e.g, \cite{brown}). In this regards, note that
	the existence results in Theorems \ref{thmCC2}, \ref{thmCC3}, \ref{thmF} for $\lambda \in [\lambda^*, \lambda^f]$ correspond to the case $\mathcal{N}^0_\lambda \neq \emptyset$. 
\end{rem}
\begin{rem}
For some other type of equations, results similar in Theorems \ref{thmCC2}, \ref{thmCC3}, that is in the case $\mathcal{N}^0_\lambda \neq \emptyset$, can also be  obtained using spectral analysis with respect to the fibering procedure   (see, for instance, \cite{ilKaye, Kaye1, Kaye2}).
\end{rem}

\section{Recursive application of the nonlinear generalized Rayleigh quotient method} \label{sec:4}

In this section, we consider a case when the fibering function $\Phi_{\lambda}(t u)$ may have more than two critical points for some $u \in W$ and $\lambda \in \mathbb{R}$.  In this case,  a direct application of the nonlinear generalized Rayleigh quotient method may not be  sufficient to find Nehari manifold extreme values. Moreover,  for such problems, it may require taking into account more than one of parameters of the problem. Below we show how this difficulty can be overcome by using,  introduced in \cite{MarcCarlIl},  recursive application of the nonlinear generalized Rayleigh quotient method.

\begin{figure}[ht]
	\centering
	\includegraphics[width=0.4\linewidth]{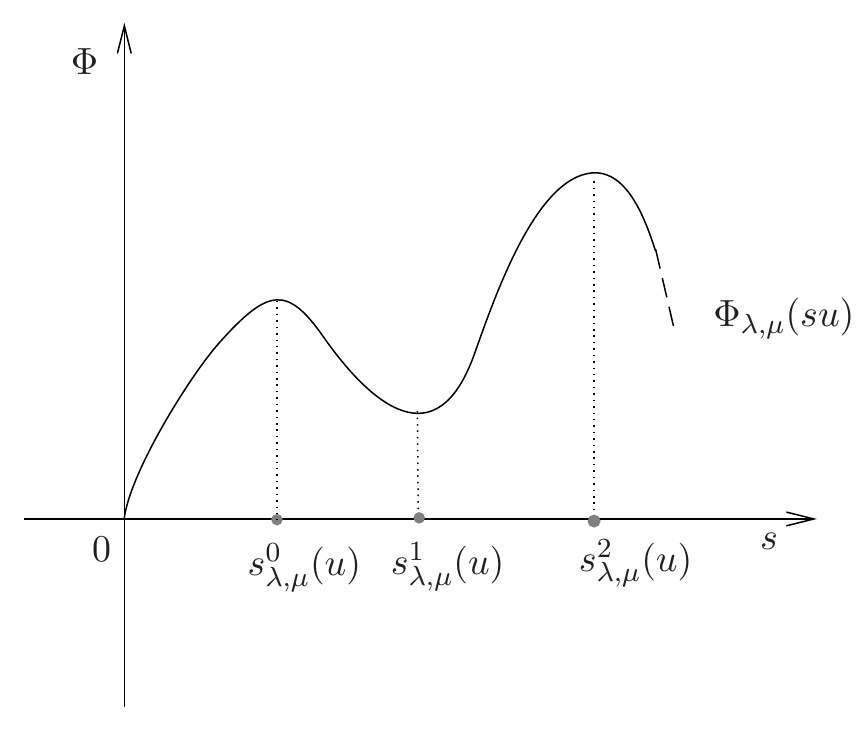}\\
	\caption{{\it Fibering function $\Phi_{\lambda,\mu}(su)$, $s\geq 0$,  $u \in W^{1,p}_0 $.}}
	\label{fig1}
\end{figure}

Let $\Omega\subset \mathbb{R}^N$ be a bounded smooth domain, $N\geq 1$. Consider the following  boundary value problem
\begin{equation}\label{p}
\left\{
\begin{aligned}
	-&\Delta_p u= |u|^{\gamma-2}u+\mu|u|^{\alpha-2}u-\lambda |u|^{q-2}u &&\mbox{in}\ \ \Omega, \\
&u=0                                   &&\mbox{on}\ \ \partial\Omega.
\end{aligned}
\right.
\end{equation}
Here $1<q<\alpha<p<\gamma<p^*$ and $\lambda,\mu \in \mathbb{R}$. By a weak solution of \eqref{p} we mean a critical point $u \in W^{1,p}_0$ of the energy functional
\begin{equation*}
\Phi_{\lambda,\mu}(u)=\frac{1}{p}\int |\nabla u|^p\,dx+\frac{\lambda}{q}\int |u|^q\,dx-\frac{\mu}{\alpha}\int |u|^\alpha-\frac{1}{\gamma}\int  |u|^\gamma\,dx.
\end{equation*}
It easily seen that for $u\in  W^{1,p}_0\setminus 0$, the fibering function $\Phi_{\lambda,\mu}(su)$  may have at most three nonzero critical points 
$$0< s^0_{\lambda,\mu}(u)\leq  s^1_{\lambda,\mu}(u) \leq s^2_{\lambda,\mu}(u)<\infty
$$ 
such that (see Figure \ref{fig1})
$$
\Phi''_{\lambda,\mu}(s^0_{\lambda,\mu}(u) u)\leq 0,~~~ \Phi''_{\lambda,\mu}(s^1_{\lambda,\mu}(u) u)\geq 0 ,~~~ \Phi''_{\lambda,\mu}(s^2_{\lambda,\mu}(u) u)\leq 0.
$$
Observe  that the case $s^j_{\lambda,\mu}(u)=s^k_{\lambda,\mu}(u)$ may occur for some $j,k\in \{0,1,2\}$, $j\neq k$,  and so $\Phi''_{\lambda,\mu}(s^j(u) u)=0$ would be true. An additional difficulty, in application of the Nehari manifold method, is that the critical points  $ s^0_{\lambda,\mu} (u)$, $ s^2_{\lambda,\mu}(u)$ satisfy the same condition $\Phi''_{\lambda,\mu}(s^0_{\lambda,\mu}(u) u)\leq 0$, $\Phi''_{\lambda,\mu}(s^2_{\lambda,\mu}(u) u)\leq 0$. 


To apply the Nehari manifold method, we need to find the values of $\lambda$, $\mu$ for which the strong inequalities $0< s^0_{\lambda,\mu}(u)< s^1_{\lambda,\mu}(u) <s^2_{\lambda,\mu}(u)$ hold.   
We solve this problem  by recursively application of the nonlinear generalized Rayleigh quotient method so that the simplest problem of zero codimension of degeneracies will be obtained at the final step of the recursion.

In the first step of the recursive procedure we consider
\begin{align}
	&\mathcal{R}^n_\lambda(u)=\frac{\int|\nabla u|^p\,dx+\lambda\int |u|^q-\int |u|^\gamma\,dx}{\int |u|^\alpha\,dx}, ~~u\in  W^{1,p}_0\setminus 0, ~\lambda\in \mathbb{R}. \label{RaylQS}
	\end{align} 
Notice that for $u\in  W^{1,p}_0\setminus 0$, $\mathcal{R}^n_\lambda(u)=\mu$ $\Leftrightarrow$ $\Phi'_{\lambda,\mu}(u)=0$. Furthermore, it readily cheek that $\mathcal{R}^n_\lambda(u)=\mu$, $\mathcal{R}^n_\lambda(u)<0\,(>0)\,(=0)$ $\Leftrightarrow$ $\Phi'_{\lambda,\mu}(u)=0$, $\Phi''_{\lambda,\mu}(u)<0\,(>0)\,(=0)$.

\begin{figure}[!ht]
\begin{minipage}[h]{0.49\linewidth}
\center{\includegraphics[scale=0.7]{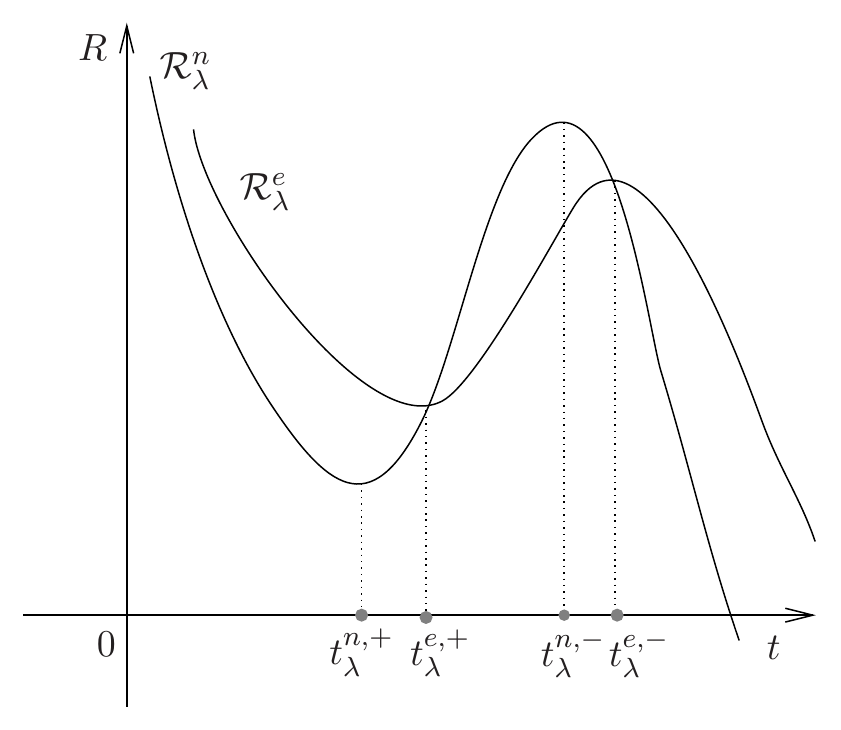}}
\caption{The functions $\mathcal{R}^e_\lambda(tu)$, $\mathcal{R}^n_\lambda(tu)$}
\label{fig2}
\end{minipage}
\hfill
\begin{minipage}[h]{0.49\linewidth}
\center{\includegraphics[scale=0.7]{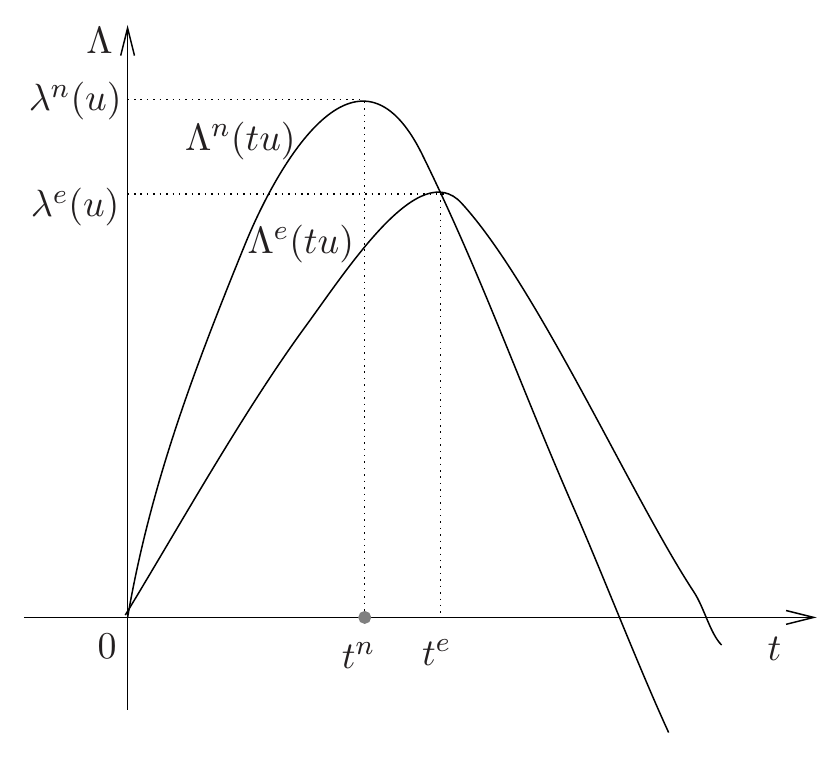}}
\caption{The functions $\Lambda^n(tu)$, $\Lambda^e(tu)$}
\label{fig3}
\end{minipage}
\end{figure}

Simple analysis shows that for any $u \in W^{1,p}_0\setminus 0$, the fibering function $ \mathcal{R}^n_\lambda(tu)$ may have at most two non-zero fibering critical points such that $0<t_\lambda^{n,+}(u)\leq t_\lambda^{n,-}(u)<+\infty$, where $t_\lambda^{n,+}(u)$ is a local minimum, $t_\lambda^{n,-}(u)$ is a local maximum point of $ \mathcal{R}^n_\lambda(tu)$, and 
\begin{equation}\label{tsts}
	0<s^0_{\lambda,\mu}(u)\leq t_\lambda^{n,+}(u) \leq s^1_{\lambda,\mu}(u)\leq t_\lambda^{n,+}(u) \leq s^2_{\lambda,\mu}(u)<\infty.
\end{equation}
Thus, $ (\mathcal{R}^n_\lambda)'(t_\lambda^{n,\pm}(u)u)=0$ and $\mathcal{R}^n_\lambda(s^k_{\lambda,\mu}(u)u)=\mu$, $k=0,1,2$. 
(see Figure \ref{fig2}).

In the second step of the recursive procedure, we apply the nonlinear generalized Rayleigh quotient method to the functional $\mathcal{R}^n_\lambda$  with respect to the parameter  $\lambda$, i.e., we consider 
\begin{equation}\label{110}
\Lambda^n(u):= \frac{(p-\alpha) \int|\nabla u|^p\,dx-(\gamma-\alpha)\int |u|^\gamma\,dx}{(\alpha-q)\int |u|^q\,dx}, ~~u \in W^{1,p}_0\setminus 0.
\end{equation}
Notice that for any $u \in W^{1,p}_0\setminus 0$, $(\mathcal{R}^n_\lambda)'(tu)=0$
$\Leftrightarrow$ $\Lambda^n(tu)=\lambda$.
Observe that the only solution of
$
\displaystyle{\frac{d}{dt}\Lambda^n(tu)=0}
$
is the   global maximum point of the function $\Lambda^n(tu)$ (see Figure \ref{fig3}) defined by  
\begin{equation*}
t^n(u):= \left(C_n\frac{\int|\nabla u|^p\,dx}{\int |u|^\gamma\,dx}\right)^{1/(\gamma-p)},~~\forall u\in  W^{1,p}_0\setminus 0,
\end{equation*}
where 
$$
C_n=\frac{(p-\alpha)(p-q)}{(\gamma-\alpha)(\gamma-q)}.
$$ 
Consider the corresponding  \textit{nonlinear generalized Rayleigh $\lambda$-quotient}   
\begin{equation}\label{const0}
\lambda^n(u):=\Lambda^n(t^n(u)u)=c^n_{q,\gamma}\frac{(\int|\nabla u|^p\,dx)^{\frac{\gamma-q}{\gamma-p}}}{(\int |u|^q\,dx) (\int |u|^\gamma\,dx)^{\frac{p-q}{\gamma-p}}},
\end{equation}
which has the following \textit{nonlinear generalized Rayleigh $\lambda$-extremal value}   
\begin{align}\label{lambExtrem}
&\lambda^n=\inf_{u\in   W^{1,p}_0\setminus 0}\sup_{t>0}\Lambda^n(tu),
=	 c^n_{q,\gamma}\frac{(\int|\nabla u|^p\,dx)^{\frac{\gamma-q}{\gamma-p}}}{(\int |u|^q\,dx) (\int |u|^\gamma\,dx)^{\frac{p-q}{\gamma-p}}},
\end{align}
where 
$$
	c^n_{q,\gamma}= \frac{(p-\alpha)^{\frac{\gamma-q}{\gamma-p}}(p-q)^{\frac{p-q}{\gamma-q}}(\gamma-p)}{(\alpha-q)(\gamma-\alpha)^{\frac{p-q}{\gamma-p}}(\gamma-q)^{\frac{\gamma-q}{\gamma-p}}}.
$$
It is not hard to show the following (see \cite{MarcCarlIl})
\begin{prop}\label{propRN2} 
	For any $\lambda\in (0,\lambda^n)$ and $u\in W^{1,p}_0\setminus 0$, the function $\mathcal{R}^n_\lambda(tu)$ has precisely two distinct critical points such that $0<t_\lambda^{n,+}(u)<t_\lambda^{n,-}(u)$. 
	Moreover,  $(\mathcal{R}^n_\lambda)''(t_\lambda^{n,+}(u)u)>0$ and  $(\mathcal{R}^n_\lambda)''(t_\lambda^{n,-}(u)u)<0$.
 \end{prop}
Observe that this and \eqref{tsts} imply that $0<s^0_{\lambda,\mu}(u)<s^2_{\lambda,\mu}(u)<\infty$ for any $\lambda\in (0,\lambda^n)$. However, we  may still have $s^0_{\lambda,\mu}(u)=s^1_{\lambda,\mu}(u)$ or $s^1_{\lambda,\mu}(u)=s^2_{\lambda,\mu}(u)$, and thus $\Phi''_{\lambda,\mu}(s^k_{\lambda,\mu}(u) u)= 0$, $k=0,1,2$. 

Notice that for $\lambda\in (0,\lambda^n)$ we are able to introduce the following \textit{nonlinear generalized Rayleigh $\mu$-quotients}
\begin{equation*}\label{mumuN}
\mu^{n,+}_\lambda(u):=\mathcal{R}^n_\lambda(t_\lambda^{n,+}(u)u),~~ \mu^{n,-}_\lambda(u):=\mathcal{R}^n_\lambda(t_\lambda^{n,-}(u)u), ~~u\in  W^{1,p}_0\setminus 0, 
\end{equation*}
resulting in the following  \textit{nonlinear generalized Rayleigh  $\mu$-extremal values} 
	\begin{align}
		&\mu^{n,+}_\lambda=\inf_{u\in   W^{1,p}_0\setminus 0}\mu^{n,+}_\lambda(u), \label{mupl}\\
		&\mu^{n,-}_\lambda=\inf_{u\in   W^{1,p}_0\setminus 0}\mu^{n,-}_\lambda(u). \label{muminus}
	\end{align}
\begin{prop}\label{propRN23} 
	Assume that $\lambda\in (0,\lambda^n)$ and $\mu <\mu^{n,-}_\lambda$. Then there exists $s^2_{\lambda,\mu}(u)>t_\lambda^{n,-}(u)$ such that $(\mathcal{R}^n_\lambda)'(s^2_{\lambda,\mu}(u) u)< 0$, i.e., $\Phi''_{\lambda,\mu}(s^2_{\lambda,\mu}(u) u)< 0$ for any $u \in W^{1,p}_0\setminus 0$.
 \end{prop}
\begin{proof}
Let $u \in W^{1,p}_0\setminus 0$. If 	$\mu <\mu^{n,-}_\lambda$, then $\mu <\mu^{n,-}_\lambda\leq \mu^{n,-}_\lambda(u)=\mathcal{R}^n_\lambda(t_\lambda^{n,-}(u)u)$. Hence and since the function $\mathcal{R}^n_\lambda(tu)$ monotone decreases on $(t_\lambda^{n,-}(u),+\infty)$, there exists a solution  $s^2_{\lambda,\mu}(u) \in (t_\lambda^{n,-}(u),+\infty)$ of the equation $\mathcal{R}^n_\lambda(su)=\mu$ such that $(\mathcal{R}^n_\lambda)'(s^2_{\lambda,\mu}(u) u)< 0$,  and consequently, 
$\Phi''_{\lambda,\mu}(s^2_{\lambda,\mu}(u) u)< 0$.
\end{proof}
However, this assertion is not sufficient to construct a solution corresponding to the critical point $s^2_{\lambda,\mu}(u)$ by the Nehari manifold method since $\Phi_{\lambda,\mu}(s u)$ may have another critical point $s^0_{\lambda,\mu}(u)$ which satisfies $\Phi''_{\lambda,\mu}(s^0_{\lambda,\mu}(u) u)\leq 0$.

This difficulty we overcome by using the following additional Rayleigh quotient
\begin{align*}
	&\mathcal{R}_\lambda^e(u)=\frac{\frac{1}{p}\int|\nabla u|^p\,dx+\frac{\lambda}{q}\int |u|^q\,dx-\frac{1}{\gamma}\int |u|^\gamma\,dx}{\frac{1}{\alpha}\int |u|^\alpha\,dx},~~u\in W^{1,p}_0\setminus 0, 
	\end{align*}
	which is characterized by the fact that $\mathcal{R}^e_\lambda(u)=\mu$ $\Leftrightarrow$ $\Phi_{\lambda,\mu}(u)=0$. The Rayleigh quotient $\mathcal{R}_\lambda^e(u)$ possesses similar properties to that $\mathcal{R}_\lambda^n(u)$. In particular, the fibering function $ \mathcal{R}^e_\lambda(tu)$ may have at most two non-zero fibering critical points $0<t_\lambda^{e,+}(u)\leq t_\lambda^{e,-}(u)<+\infty$ so that $t_\lambda^{e,+}(u)$ is a local minimum,  whereas $t_\lambda^{e,-}(u)$ is a local maximum point of $ \mathcal{R}^e_\lambda(tu)$. Moreover, the same conclusion as for $\Lambda^n(u)$ can be drawn for the Rayleigh quotient 
\begin{equation}\label{11E}
\Lambda^e(u):= q\frac{\frac{(p-\alpha)}{p} \int|\nabla u|^p\,dx-\frac{(\gamma-\alpha)}{\gamma}\int |u|^\gamma\,dx}{(\alpha-q)\int |u|^q\,dx}.
\end{equation}
which is characterized by the fact that   $(\mathcal{R}^e_\lambda)'(tu)=0$ $\Leftrightarrow$ $\Lambda^e(tu)=0$, for $u \in W^{1,p}_0\setminus 0$.
The unique solution of
$
\frac{d}{dt}\Lambda^e(tu)=0
$
is a  global maximum point of the function $\Lambda^e(tu)$ defined by  
\begin{equation}
\label{17}
t^e(u):= \left(C_e\frac{\| u\|^p_1}{\| u\|^\gamma_{L^\gamma}}\right)^{1/(\gamma-p)},~~\forall u\in  W^{1,p}_0\setminus 0,
\end{equation}
where 
$$
C_e=\frac{\gamma(p-\alpha)(p-q)}{p(\gamma-\alpha)(\gamma-q)}.
$$ 
Thus we have the following additional \textit{nonlinear generalized Rayleigh $\lambda$-extremal value}   
\begin{align} 
			 &\lambda^e=\inf_{u\in  W^{1,p}_0\setminus 0}\sup_{t>0}\Lambda^e(tu)=c^e_{q,\gamma}\inf_{u\in  W^{1,p}_0\setminus 0}\frac{\| u\|_1^{p\frac{\gamma-q}{\gamma-p}}}{\| u\|^q_{L^q} \| u\|_{L^\gamma}^{\gamma\frac{p-q}{\gamma-p}}},. \label{lambExtrem1}
	\end{align}	
which  makes it possible to split the extremal points of the functionals $\mathcal{R}^e_\lambda$, indeed, we have
\begin{prop}\label{propRNE22} 
	For any $\lambda\in (0,\lambda^e)$ and $u\in W^{1,p}_0\setminus 0$, the function $\mathcal{R}^e_\lambda(tu)$ has precisely two distinct critical points such that $0<t_\lambda^{e,+}(u)<t_\lambda^{e,-}(u)$. Moreover,  $(\mathcal{R}^e_\lambda)''(t_\lambda^{e,+}(u)u)>0$,  $(\mathcal{R}^e_\lambda)''(t_\lambda^{e,-}(u)u)<0$.
\end{prop} 
Thus, for $\lambda\in (0,\lambda^e)$, we have the following \textit{nonlinear generalized Rayleigh $\mu$-quotients}
\begin{equation*}\label{mumuN2}
\mu^{e,+}_\lambda(u):=\mathcal{R}^e_\lambda(t_\lambda^{e,+}(u)u),~~~~ \mu^{e,-}_\lambda(u):=\mathcal{R}^n_\lambda(t_\lambda^{e,-}(u)u), ~~u\in  W^{1,p}_0\setminus 0
\end{equation*}
with corresponding principal \textit{nonlinear generalized Rayleigh  $\mu$-extremal values} 
	\begin{align}
		&\mu^{e,+}_\lambda=\inf_{u\in   W^{1,p}_0\setminus 0}\mu^{e,+}_\lambda(u),~~\mu^{e,-}_\lambda=\inf_{u\in   W^{1,p}_0\setminus 0}\mu^{e,-}_\lambda(u).
	\end{align}
The relationships among the NG-Rayleigh quotients   is given by the following lemma

\begin{lem}\label{lemGEOM}
Assume that $1<q<\alpha<p<\gamma$, $u\in W^{1,p}_0\setminus 0$, $t>0$. 
		\begin{enumerate}
			\item[\rm{$(i)$}] $\Lambda^{e}(tu)=\Lambda^n(tu)$   $\Leftrightarrow $ $t=t^e(u)$,
			\item[\rm{$(ii)$}]    $\mathcal{R}_\lambda^e(tu)=\mathcal{R}_\lambda^n(tu)$ $\Leftrightarrow $ $t=t_\lambda^{e,+}(u)$ or $t=t_\lambda^{e,-}(u)$, $\forall \lambda \in(0,\lambda^e)$,
			\item[\rm{$(iii)$}]   $t_\lambda^{n,+}(u)<t_\lambda^{e,+}(u)<t^e(u)<t_\lambda^{n,-}(u)<t_\lambda^{e,-}(u)$, $\forall \lambda \in(0,\lambda^e)$.
\end{enumerate}
	\end{lem}
	\begin{proof} The equality
				$\Lambda^{e}(tu)=\Lambda^n(tu)$  is equivalent to 
				$$
				 t^{p-q}\| u\|^p_1-\frac{(\gamma-\alpha)}{(p-\alpha)}t^{\gamma-q}\|u\|^\gamma_{L^\gamma}=q\left(t^{p-q}\frac{1}{p} \| u\|^p_1-t^{\gamma-q}\frac{(\gamma-\alpha)}{\gamma(p-\alpha)}\|u\|^\gamma_{L^\gamma}\right).
				$$
Hence,
				$$
				0=\frac{(p-q)(p-\alpha)}{p} t^{1-q}\| u\|^p_1-\frac{(\gamma-q)(\gamma-\alpha)}{\gamma}t^{\gamma-q-1}\|u\|^\gamma_{L^\gamma}=(\Lambda^{e}(tu))',
				$$
which implies \rm{$(i)$}.
		
Observe, $\mathcal{R}_\lambda ^e(tu)=\mathcal{R}_\lambda^n(tu)$ for $t>0$ if and only if 
				$$t^{p-\alpha}\| u\|^p_1+\lambda t^{q-\alpha}\|u\|^q_{L^q}-t^{\gamma-\alpha}\|u\|^\gamma_{L^\gamma}=\frac{\alpha t^{p-\alpha}}{p}\| u\|^p_1+\frac{\lambda\alpha t^{q-\alpha}}{q}-\frac{\alpha t^{\gamma-\alpha}}{\gamma}\|u\|^\gamma_{L^\gamma},$$
which is equivalent to
						\begin{align*}
					0&=\frac{(p-\alpha)}{p}t^{p-\alpha}\| u\|^p_1-\frac{\gamma-\alpha}{\gamma}t^{\gamma-\alpha}\|u\|^\gamma_{L^\gamma}-\frac{\lambda(\alpha-q)}{q}t^{q-\alpha}\|u\|^q_{L^q}\\
					&=\frac{(\alpha-q)\|u\|^q_{L^q}t^{q-\alpha}}{q}\left(\Lambda_{e}(tu)-\lambda\right).
				\end{align*}
Since $(\mathcal{R}^e_\lambda)'(tu)=0$ $\Leftrightarrow$ $\Lambda^e(tu)=0$, we get \rm{$(ii)$}.	
The proof of  \rm{$(iii)$} follows from \rm{$(i)$},\rm{$(ii)$}.	
\end{proof} 
Lemma \ref{lemGEOM} can be illustrated as follows by Figures \ref{fig2}, \ref{fig3}.
\begin{lem} \label{TH0}
Assume that $1<q<\alpha<p<\gamma<p^*$. Then
	\begin{description}
		\item[($1^o$)] 
		 $ 0<\lambda^e<\lambda^n<+\infty$,
		\item[($2^o$)]   $0<  \mu^{n,+}_\lambda< \mu^{e,+}_\lambda<\mu^{e,-}_\lambda<\mu^{n,-}_\lambda<+\infty$ for  $\lambda \in  (0,\lambda^e)$,
		\item[($3^o$)]  $0< \mu^{n,+}_\lambda<\mu^{n,-}_\lambda<+\infty$ for $\lambda \in  (0,\lambda^n)$.  
\end{description}
\end{lem}
For the proof we refer the reader to  \cite{MarcCarlIl}.

\begin{thm}\label{thm11}
Assume that $1<q<\alpha<p<\gamma<p^*$. Then 

\par\noindent	
$(1^o)$ For any $\lambda \in (0,\lambda^e)$, $\mu\in (\mu^{e,+}_\lambda,\mu^{n,-}_\lambda)$, problem \eqref{p}
possesses a weak positive solution $u^1_{\lambda, \mu}\in   C^{1,\kappa}(\overline{\Omega})$,  $\kappa \in (0, 1)$ such that 
\begin{equation}\label{QP1}
	\Phi_{\lambda,\mu}(u^1_{\lambda, \mu})<0,~~\Phi''_{\lambda,\mu
	}(u^1_{\lambda, \mu})>0.
\end{equation}
Moreover, 	$u^1_{\lambda, \mu}$ is a ground state of \eqref{p}.

\par\noindent	
$(2^o)$	For any $\lambda \in (0,\lambda^n)$, $\mu \in (\mu^{e,-}_\lambda, \mu^{n,-}_\lambda)$,  problem \eqref{p} has a weak positive solution $u^2_{\lambda, \mu}\in   C^{1,\kappa}(\overline{\Omega})$,  $\kappa \in (0, 1)$ such that
		\begin{equation}\label{QP}
		\Phi_{\lambda,\mu}(u^2_{\lambda, \mu})<0, ~~\Phi''_{\lambda,\mu}(u^2_{\lambda, \mu})<0.
	\end{equation}
	
\end{thm}
\begin{proof}
Consider the Nehari manifold 
\begin{equation*}\label{eq:NehariR}
		\mathcal{N}_{\lambda,\mu}=\{u \in W^{1,p}_0\setminus 0: ~\mathcal{R}^n_\lambda(u)=\mu\}.
\end{equation*}
Observe that 	$\Phi_{\lambda,\mu}$ is coercive on 
		$\mathcal{N}_{\lambda,\mu}$, $\forall \lambda >0$, $\forall\mu \in \mathbb{R}$. Indeed, by the Sobolev inequality,
\begin{align*}\label{Phi-bounded}
							\Phi_{\lambda,\mu}(u)\geq \frac{\gamma-p}{p\gamma}\|u\|_1^p-\mu C\|u\|_1^\alpha,~~\forall u\in \mathcal{N}_{\lambda,\mu},
							\end{align*}
where $C<+\infty$  does not depend on $u \in \mathcal{N}_{\lambda,\mu}$. Since $\alpha<p$, this implies $\Phi_{\lambda,\mu}(u) \to +\infty$ if  $\|u\|_1 \to +\infty$ and $u \in \mathcal{N}_{\lambda,\mu}$.

\noindent $(1^o)$ Denote 
$$
\mathcal{R}\mathcal{N}_{\lambda,\mu}^1:=\{u\in W^{1,p}_0\setminus 0:~ (\mathcal{R}_\lambda^{n})'(u)>0,~ u \in \mathcal{N}_{\lambda,\mu}\}.
$$ 
Consider
\begin{equation} \label{GminUnst0}
	\hat{\Phi}^1_{\lambda,\mu}=\min \{\Phi_{\lambda,\mu}(u):  ~ u \in \mathcal{R}\mathcal{N}_{\lambda,\mu}^1\}.
\end{equation}
 Observe $\hat{\Phi}^1_{\lambda,\mu}<0$. Indeed,  $\mu<\mu^{n,-}_\lambda<\mu^{n,-}_\lambda(u)$, $\forall u \in W^{1,p}_0\setminus 0$.  Hence, for $\mu \in (\mu^{e,+}_\lambda, \mu^{n,-}_\lambda)$, there exists $\tilde{u} \in W^{1,p}_0\setminus 0$ such that $\mu^{e,+}_\lambda(\tilde{u})<\mu<\mu^{n,-}_\lambda(\tilde{u})$.  Lemma \ref{lemGEOM} implies that  $s^1_{\lambda, \mu}(\tilde{u}) \in  (t^{e,+}_\lambda(\tilde{u}), t^{e,-}_\lambda(\tilde{u}))$ and  $\mathcal{R}_\lambda^n(t\tilde{u})>\mathcal{R}_\lambda^e(t\tilde{u})$ for $t \in (t^{e,+}_\lambda(\tilde{u}), t^{e,-}_\lambda(\tilde{u}))$.  Hence, $\mu=\mathcal{R}_\lambda^n(s^1_{\lambda,\mu}(\tilde{u})\tilde{u})>\mathcal{R}_\lambda^e(s^1_{\lambda,\mu}(\tilde{u})\tilde{u})$, and thus   $\hat{\Phi}^1_{\lambda,\mu}\leq \Phi_{\lambda, \mu}(s^1_{\lambda,\mu}(\tilde{u})\tilde{u}) <0$. In addition, we have got that $\mathcal{R}\mathcal{N}_{\lambda,\mu}^1\neq \emptyset$, $\forall \lambda \in (0,\lambda^e)$, $\forall \mu \in (\mu^{e,+}_\lambda, \mu^{n,-}_\lambda)$.

	Let $u_m$ be a minimizing sequence of \eqref{GminUnst0}, i.e.,
	$$
	\Phi_{\lambda,\mu}(u_m) \to \hat{\Phi}^1_{\lambda,\mu},~~\Phi_{\lambda,\mu}'(u_m)=0,~ (\mathcal{R}_\lambda^{n})'(u_m)>0,~~m=1,\ldots .
	$$
By the coercivity of $\Phi_{\lambda,\mu}$  on $\mathcal{N}_{\lambda,\mu}$, the sequence $(u_m)$ is bounded in $W^{1,p}_0$ and thus,  up to a
subsequence, 
$$
u_m \to \bar{u}~~\mbox{strongly in}~~L^r,~~\mbox{and weakly in} ~W^{1,p}_0,
$$
for some $\bar{u} \in W^{1,p}_0$, where $r \in(1,p^*)$. Observe, $\Phi_{\lambda,\mu}(\bar{u})\leq \liminf_{m\to \infty}\Phi_{\lambda, \mu}(u_m)=\hat{\Phi}^1_{\lambda,\mu}<0$, and thus $\bar{u}\neq 0$.

Notice that  $\lambda<\lambda^e<\lambda^n<\lambda^e(\bar{u})<\lambda^n(\bar{u})$. Consequently, there exist nonlinear generalized Rayleigh $\mu$-quotients $\mu^{e,+}_\lambda(\bar{u}),\mu^{n,-}_\lambda(\bar{u})$ such that $\mu^{e,+}_\lambda(\bar{u})<\mu^{n,-}_\lambda(\bar{u})$. Since $\Phi_{\lambda,\mu}(\bar{u})<0$, we have $\mathcal{R}^e_\lambda(\bar{u})<\mu$, which implies  $\mu^{e,+}_\lambda(\bar{u})\leq \mathcal{R}^e_\lambda(\bar{u})<\mu$. Thus  $\mu^{e,+}_\lambda(\bar{u})<\mu<\mu^{n,-}_\lambda(\bar{u})$, and 
there exists $s^1_{\lambda, \mu}(\bar{u})\in (s^0_{\lambda, \mu}(\bar{u}), s^2_{\lambda, \mu}(\bar{u}))$ such that  
\begin{equation*}
\begin{aligned}
	&\mu=\mathcal{R}^n_\lambda(s^1_{\lambda, \mu}(\bar{u})\bar{u})\leq \liminf_{m\to \infty}\mathcal{R}^n_\lambda(s^1_{\lambda, \mu}(\bar{u})u_m),\\
		&0<(\mathcal{R}^n_\lambda)'(s^1_{\lambda, \mu}(\bar{u})\bar{u})\leq  \liminf_{m\to \infty}(\mathcal{R}^n_\lambda)'(s^1_{\lambda, \mu}(\bar{u})u_m).
			\end{aligned}
\end{equation*}
This means that $1=s^1_{\lambda,\mu}(u_m)\leq s^1_{\lambda, \mu}(\bar{u}) <s^2_{\lambda, \mu}(u_m)$, $m=1,\ldots $.  Since $\Phi'_{\lambda, \mu}(t u_m)<0$ for any $t \in (1,s^1_{\lambda, \mu}(\bar{u}))$,
$$
\Phi_{\lambda,\mu}(s^1_{\lambda, \mu}(\bar{u})\bar{u})\leq \Phi_{\lambda,\mu}(\bar{u})\leq \liminf_{m\to \infty}\Phi_{\lambda, \mu}(u_m)=\hat{\Phi}^1_{\lambda,\mu},
$$
which yields that  $\bar{u}$ is a minimizer of  \eqref{GminUnst0}. In view of that $\mu<\mu^{n,-}_\lambda<\mu^{n,-}_\lambda(\bar{u})$, 
we obtain that $u^1_{\lambda,\mu}:=\bar{u}^1$ weakly satisfies  \eqref{p}.   A trivial verification shows that $u^1_{\lambda,\mu}$ is a ground state. The rest of the proof is similar to that in Theorem \ref{cor444}.

\noindent
  $(2^o)$  
Consider
\begin{equation} \label{GminUnst1}
	\hat{\Phi}^2_{\lambda,\mu}=\min \{\Phi_{\lambda,\mu}(u)~:~ u \in \mathcal{R}\mathcal{N}_{\lambda,\mu}^{2}\},
\end{equation}
where 
$$
\mathcal{R}\mathcal{N}_{\lambda,\mu}^{2}=\{u\in W^{1,p}_0\setminus 0: \mathcal{R}_{\lambda}^n(u)=\mu,~ (\mathcal{R}^n_\lambda)'(u)<0,~\mathcal{R}_\lambda^e(u)<\mu\}.
$$
The assumption $\mu \in (\mu^{e,-}_\lambda, \mu^{n,-}_\lambda)$ implies that there exists $\tilde{u} \in W^{1,p}_0\setminus 0$ such that $\mu^{e,-}_\lambda(\tilde{u})<\mu<\mu^{n,-}_\lambda(\tilde{u})$.  Hence, $\mu>\mu^{e,-}_\lambda(\tilde{u})>\mathcal{R}_{\lambda}^e(s^2(\tilde{u})\tilde{u})$, and thus,   $s^2(\tilde{u})\tilde{u} \in  \mathcal{R}\mathcal{N}_{\lambda,\mu}^{2} \neq \emptyset$, $\forall \lambda \in (0,\lambda^n)$, $\forall \mu \in (\mu^{e,-}_\lambda, \mu^{n,-}_\lambda)$.
Furthermore, the inequality $\mathcal{R}_\lambda^e(u)<\mu$, $u \in \mathcal{R}\mathcal{N}_{\lambda,\mu}^{2}$ implies $\Phi_{\lambda,\mu}(u)<0$ and consequently $\hat{\Phi}^2_{\lambda,\mu}<0$.

Let $(u_m)$ be a minimizing sequence of \eqref{GminUnst1}. Similar to the proof of $(1^o)$ one deduces that there exists a subsequence, which we again denote by $(u_m)$, and a nonzero limit point $\bar{u}^2$ such that
					$$
					u_m \to \bar{u}^2~~\mbox{strongly in}~L^r,~r \in (1,p^*), ~~\mbox{and weakly in} ~W^{1,p}_0.
					$$
As above, it follows that $\lambda<\lambda^n(\bar{u})$ and  $\mu<\mu^{n,-}_\lambda(\bar{u}^2)$. Consequently, there exist $s^2_{\lambda, \mu}(\bar{u})> 0$ and $t^{n,\pm}_{\lambda}(\bar{u}^2)>0$ such that  $t^{n,+}_{\lambda}(\bar{u}^2)<t^{n,-}_{\lambda}(\bar{u}^2)$, and
					\begin{equation*}\label{ineqPHI}
\begin{aligned}
	&\mu=\mathcal{R}^n_\lambda(s^2_{\lambda, \mu}(\bar{u}^2)\bar{u}^2)\leq \liminf_{m\to \infty}\mathcal{R}^n_\lambda(s^2_{\lambda, \mu}(\bar{u}^2)u_m),\\
		&0=(\mathcal{R}^n_\lambda)'(t^{n,\pm}_{\lambda}(\bar{u}^2)\bar{u}^2)\leq \liminf_{m\to \infty}(\mathcal{R}^n_\lambda)'(t^{n,\pm}_{\lambda}(\bar{u}^2)u_m),
			\end{aligned}
\end{equation*}
where the first inequality  implies that $s^2_{\lambda, \mu}(\bar{u}^2)\in (0,s^0_{\lambda,\mu}(u_m))\cup (s^1_{\lambda,\mu}(u_m),s^2_{\lambda,\mu}(u_m)=1)$, whereas from the second one we get $t^{n,-}_{\lambda}(\bar{u}^2)\in (t^{n,+}_{\lambda}(u_m),t^{n,-}_{\lambda}(u_m)) $,  $m=1,\ldots$. Hence, $s^0_{\lambda,\mu}(u_m)<t^{n,+}_{\lambda}(u_m)<t^{n,-}_{\lambda}(\bar{u}^2)<s^2_{\lambda, \mu}(\bar{u}^2)$, and thus $s^2_{\lambda, \mu}(\bar{u}^2)\in  (s^1_{\lambda,\mu}(u_m),s^2_{\lambda,\mu}(u_m))$. Since $\Phi'_{\lambda, \mu}(tu_m)>0$ for any $s^1_{\lambda, \mu}(u_m)<t<s^2_{\lambda, \mu}(u_m)$,
\begin{equation}\label{FGG}
	\Phi_{\lambda,\mu}(s^2_{\lambda, \mu}(\bar{u}^2)\bar{u}^2)\leq\liminf_{m\to \infty}\Phi_{\lambda,\mu}(s^2_{\lambda, \mu}(\bar{u}^2)u_m)\leq \liminf_{m\to \infty}\Phi_{\lambda,\mu}(s^2_{\lambda, \mu}(u_m)u_m)=\hat{\Phi}^{2}_{\lambda,\mu}<0.
\end{equation}
The inequality $\Phi_{\lambda,\mu}(s^2_{\lambda, \mu}(\bar{u}^2)\bar{u}^2)<0$ implies that $\mathcal{R}_{\lambda}^e(s_2(\bar{u}^2)\bar{u}^2)<\mu$. Thus,  $s_2(\bar{u}^2)\bar{u}^2 \in \mathcal{R}\mathcal{N}_{\lambda,\mu}^{2}$, and  $s_2(\bar{u}^2)\bar{u}^2$ is a minimizer of  \eqref{GminUnst1}. By \eqref{FGG}, this implies that  $u_m \to \bar{u}^2$ strongly in $W^{1,p}_0$  and $s^2_{\lambda, \mu}(\bar{u}^2)=1$. Consequently, $u^2_{\lambda,\mu}:= \bar{u}^2$ is a minimizer  of \eqref{GminUnst1}. The rest of the proof runs as before. 
	
\end{proof}
\begin{cor}
	Assume that $1<q<\alpha<p<\gamma<p^*$. Then for any $\lambda \in (0,\lambda^e)$, $\mu\in (\mu^{e,-}_\lambda,\mu^{n,-}_\lambda)$, problem \eqref{p}
possesses two distinct weak positive solutions $u^1_{\lambda, \mu}, u^2_{\lambda, \mu}\in   C^{1,\kappa}(\overline{\Omega})$,  $\kappa \in (0, 1)$ such that
\begin{equation*}
	\Phi_{\lambda,\mu}(u^1_{\lambda, \mu})<0,~~\Phi''_{\lambda,\mu
	}(u^1_{\lambda, \mu})>0,~~\Phi_{\lambda,\mu}(u^2_{\lambda, \mu})<0, ~~\Phi''_{\lambda,\mu}(u^2_{\lambda, \mu})<0.
\end{equation*}
Moreover, 	$u^1_{\lambda, \mu}$ is a ground state of \eqref{p}.
\end{cor}

\begin{figure}[ht]
	\centering
	\includegraphics[width=0.5\linewidth]{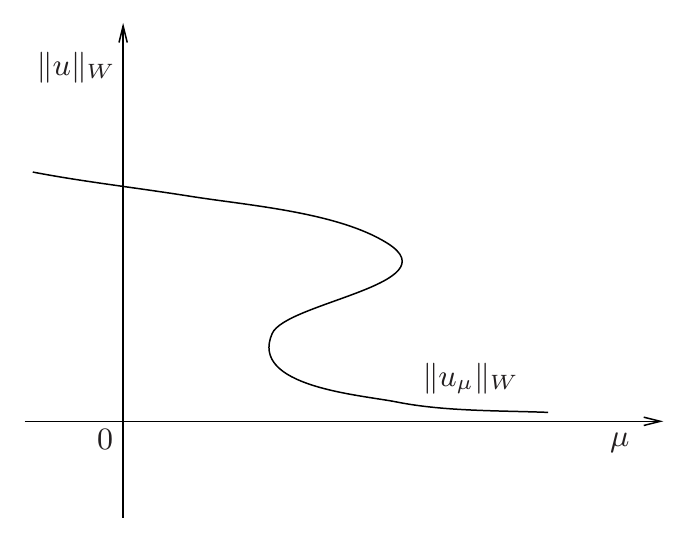}\\
	\caption{{\it The branches of solutions to \eqref{p} with fix $\lambda\in (0,\lambda^e)$.}}
	\label{fig4}
\end{figure}

\begin{rem}
We anticipate that using the mountain pass theorem one can show that \eqref{p} has a third positive solution. This leads us to the conjecture that the branch of solutions to \eqref{p} has of S-shape type bifurcation curve, see Figure \ref{fig4}.  
\end{rem}
\begin{rem}
Results concerning the existence of two or more nonnegative solution have been obtained  by using variational approaches for some another type of equations, as well. For instance, Lubushev in \cite{Lubyshev}, Bobkov, Dr\'abek,  Hern\'andez in \cite{BobDrHer}  using   global and local minimizations, and  fibering method and mountain pass theorem as well, obtained the existence of three positive solution for equation \eqref{p} in the case $1<q<p<\alpha<\gamma<p^*$, for some $\lambda, \mu \in (-\infty,0)$. However, the geometry of the fibering function $\Phi_{\lambda,\mu}(tu)$ in this case is different than that  have been considered above.  This fact prevents a direct application of the result of \cite{BobDrHer,Lubyshev} to our problem.
\end{rem}
\begin{rem}
We wonder if one can find solution $u^2_{\lambda, \mu}$ as in  Theorem \ref{thm11} without using the NG-Rayleigh quotients $\mathcal{R}^n_\lambda(u)$, $\mathcal{R}^e_\lambda(u)$. Note that the function $\Phi'''_{\lambda,\mu}(su)$ can change sign twice on $(0, +\infty)$ if $\gamma<\alpha +1$ (see e.g. \cite{MarcCarlIl}), that can cause difficulties in the splitting the critical points of $\Phi_{\lambda,\mu}(u)$ if the NG-Rayleigh quotients are not used.
\end{rem}

\section{Prescribed energy solution of the zero frequency problem } \label{sec:5}
In this section, we use an energy level Rayleigh quotient to prove an existence of solution of the following zero-frequency problem \cite{ilStab} (which is also called  the "zero mass" case problem (cf. \cite{BerestL}):
\begin{equation}\label{1S}
	-\Delta u -\mu|u|^{p-2}u+|u|^{q-2}u =0, ~~x \in \mathbb{R}^N,   
\end{equation}
where $p, q \in (2,2^*)$, $N\geq 3$ and $\mu\in \mathbb{R}$. The problem  has a variational structure and the associated energy functional
$$
E_\mu(u)= \int \left(\frac{1}{2} |\nabla u|^2 - \frac{\mu}{p}| u|^p +\frac{1}{q}| u|^q\right)dx
$$
is defined in 
$
\mathcal{D}:=\mathcal{D}^{1,2}(\mathbb{R}^N)\cap L^q(\mathbb{R}^N)\cap L^p(\mathbb{R}^N), 
$
where
$$
\mathcal{D}^{1,2}:=\mathcal{D}^{1,2}(\mathbb{R}^N) := \{u \in L^{2^*}(\mathbb{R}^N): \nabla u \in L^{2} (\mathbb{R}^N)\}
$$
the space with  norm $\|u\|_{\mathcal{D}^{1,2}}:=\int |\nabla u |^2 \, dx$ (cf. \cite{BerestL, willem}).

 For $u \in \mathcal{D}^{1,2}$, we denote $u_\sigma:=u(x/\sigma)$, $x \in \mathbb{R}^N$,  $\sigma>0$, and 
$$
 T(u):=\int|\nabla u|^{2}\, dx,~~A(u):= \int |u|^{p}\, dx,~~B(u):=\int |u|^{q}\, dx.
$$

For $E\geq 0$, consider the energy level Rayleigh quotient
\begin{align}
	R^E(u):= \frac{\frac{1}{2}T(u)+\frac{1}{q}B(u)-E}{\frac{1}{p}A(u)}.
\end{align}

Let $u \in H^1\setminus 0$, $E>0$, $\sigma>0$, consider
\begin{equation*}
		R^E(u_\sigma)=\frac{\sigma^{-2}\frac{1}{2}T(u)+\frac{1}{q}B(u)-\sigma^{-N}E}{\frac{1}{p}A(u)}.
	\end{equation*}
Then
$$
\frac{d}{d\sigma} (R^E)(u_\sigma )=0 ~\Leftrightarrow~~\sigma=\sigma^E(u):=\left(\frac{NE}{T(u)}\right)^\frac{1}{N-2}.
$$
Hence,  we are able to introduce the following nonlinear generalized Rayleigh quotient 
\begin{equation}\label{LCN}
	M^E(u ):=R^E(u_{\sigma^E(u)} )=\frac{p}{A(u)}\left(\frac{c_N^E}{2}T^\frac{N}{(N-2)}(u) +\frac{1}{q}B(u)\right),
\end{equation}
where
\begin{equation*}\label{CN}
	c_N^E=\frac{(N-2)}{N^\frac{N}{(N-2)}E^\frac{2}{(N-2)}}.
\end{equation*}
Observe that $M^E(u )$ is a $0$-homogeneous functional with respect to the scale change $\sigma \mapsto u_\sigma$. Furthermore, $DM^E(u )=0$, $M^E(u )=\mu$, $\sigma^E(u)=1$ $\Leftrightarrow$ $DE_\mu(u)=0$, $E_\mu(u)=E$

For every $u \in\mathcal{D}\setminus 0$,  consider the fibering function
\begin{equation*}\label{fiberMu}
	M^E(tu):=\frac{\frac{c_N^E}{2}t^{2^*-p}T^\frac{N}{(N-2)}(u)+\frac{1}{q}t^{q-p}B(u)}{\frac{1}{p}A(u)}, ~~t>0.
\end{equation*}
Notice that  if $2<p<q<2^*$, then $(M^E)'(su)\neq 0$, $\forall t>0$, $\forall u \in \mathcal{D}$. Hence
\begin{cor}\label{NessC}
If $2<p<q<2^*$, then \eqref{1S} has no weak solution in 	$\mathcal{D}\setminus 0$ with any energy level $E\geq 0$.
\end{cor}

Assume that $2<q<p<2^*$. It is easily seen that  the function $t\mapsto  M^E(tu)$  attains its global minimum at the unique point 
$$
t^E(u)=\left(\frac{B(u)}{ c_{p,q,N,E}T^\frac{N}{(N-2)}(u)}\right)^{1/(2^*-q)},
$$
where $ c_{p,q,N,E}=c_N^E q(2^*-p)/2(p-q)$. 
Introduce
\begin{equation}\label{MSu}
	\mu^E(u):=M^E(t^E(u)u)=\min_{t\geq 0}M^E(tu)= C_{p,q,N,E}\frac{B^\frac{2^*-p}{2^*-q}(u)T^{\frac{2^*(p-q)}{2(2^*-q)}}(u)}{A(u)},
\end{equation}
where 
\begin{equation*}\label{Cpq}
	 C_{p,q,N,E}= \frac{c(p,q,N)}{E^\frac{2(p-q)}{(2^*-q)(N-2)}},
\end{equation*} 
\begin{equation}\label{cpqn}
	c(p,q,N)=\left(\frac{(N-2)}{N^\frac{N}{(N-2)}}\frac{q(2^*-p)}{2(p-q)}\right)^\frac{(p-q)}{(2^*-q)}\frac{p(2^*-q)}{q(2^*-p)}.
\end{equation}
 
\begin{cor}
	
	\begin{itemize}
		\item $D\mu^E(u )=0$, $\mu^E(u )=\mu$, $\sigma^E(u)=t^E(u)=1$  $\Leftrightarrow$ $DE_\mu(u)=0$, $E_\mu(u)=E$;
		\item  		$\mu^E(u)$ is a multi-homogeneous functional on $\mathcal{D}$, namely:  $\mu^E(u)=\mu^E(u_\sigma)=\mu^E(tu)$, $\forall \sigma>0, ~\forall t>0,~\forall u \in \mathcal{D}\setminus 0$.
\end{itemize}
\end{cor}

Denote $\beta:=\frac{2q(2^*-p)}{2^*(p-q)}$, $\rho:=\frac{2p(2^*-q)}{2^*(p-q)}$ and define 
$$
\mu(u):=\frac{\|u\|_{L^q}^{\beta}\|\nabla u\|_{L^2}^2}{\|u\|_{L^p}^{\rho}(u)}\equiv \frac{E^\frac{2(p-q)}{(2^*-q)(N-2)}}{ c(p,q,N))^{\rho/p}}\mu^E(u),~~ u \in \mathcal{D}\setminus 0.
$$
Consider  
\begin{equation}\label{MUP}
	\bar{\mu}=\inf_{ u \in \mathcal{D}\setminus 0}\mu(u).
\end{equation}
Denote
\begin{equation}
	\hat{\mu}^E:=  \frac{c(p,q,N)}{E^\frac{2(p-q)}{(2^*-q)(N-2)}}\bar{\mu}^{p/\rho},~~\forall E>0.
\end{equation}
By the Gagliardo–Nirenberg interpolation inequality 
	\begin{align}\label{GN}
		\int |u|^{p}\,dx \leq C_{gn}(\int |\nabla u|^{2}\,dx)^{\frac{2^*(p-q)}{2(2^*-q)} }&(\int |u|^{q}\,dx)^\frac{2^*-p}{2^*-q}~\Leftrightarrow \\
		&A(u)\leq C_{gn}(T(u))^{\frac{2^*(p-q)}{2(2^*-q)} }(B(u))^\frac{2^*-p}{2^*-q},\nonumber
	\end{align}
where constant $C_{gn}$ does not depend on $u \in \mathcal{D}$. Hence $0<\bar{\mu}<+\infty$.

\begin{thm}\label{ZeroLem}
Assume that $2<q<p<2^*$, $N\geq 3$. There exists a minimizer $\hat{u}^E_{\hat{\mu}^E}$ of \eqref{MUP} which is a weak solution of  \eqref{1S} with  $\mu=\hat{\mu}^E$ and prescribed energy $E_{\hat{\mu}^E}(u)=E$. Moreover,  $\hat{u}^E_{\hat{\mu}^E}>0$ in $\mathbb{R}^N$, $\hat{u}^E_{\hat{\mu}^E} \in C^{2}(\mathbb{R}^N)$.
	\end{thm}
\begin{proof}
Let $(v_i)$ be a minimizing sequence	of \eqref{MUP}, i.e.,  $\mu(v_i) \to \bar{\mu}$ as $i\to \infty$. Set $u_i=t_i\cdot(v_i)_{\sigma_i}$, $i=1,2,\ldots $,  where 
$$
t_i=(\|v_i\|_{L^q}^q/\|v_i\|_{L^p}^p)^{1/(p-q)}, ~~\sigma_i=(\|v_i\|_{L^q}^{pq}/\|v_i\|_{L^p}^{qp})^{1/N(p-q)}.
$$
Then $\|u_i\|_{L^p}=1$ and $\|u_i\|_{L^q}=1$, $i=1,\ldots $, and by the homogeneity  of $\mu(u)$, $(u_i)$ is  a minimizing sequence	of \eqref{MUP}. Since $\bar{\mu}<+\infty$, $(\|\nabla u_i\|_{L^2})$ is bounded, and thus  
$(u_i)$ is bounded in $\mathcal{D}$ and in  $W^{1,2}_{loc}(\mathbb{R}^N)$. Thus, by the Banach–Alaoglu and Sobolev embedding theorems, there exists a subsequence, still denoted by $(u_i)$, such that 
\begin{align*}
	&u_i \rightharpoonup \hat{u}^* ~~\mbox{in} ~ \mathcal{D},\\
	&u_i \to \hat{u}^*  ~\mbox{in} ~~ L^\gamma_{loc},~~1\leq \gamma<2^*,\\
	&u_i \to \hat{u}^* ~~\mbox{a.e. on} ~\mathbb{R}^N,
\end{align*}
for some $\hat{u}^* \in \mathcal{D}$. 

Let $r > 0$. Observe that 
$$
\delta:= \liminf_{i\to \infty}\sup_{y\in \mathbb{R}^N}\int_{B(y;r)} |u_i|^\gamma \,dx>0,
$$
for each $\gamma \in (1,p)$. Indeed, if this is not true, then by the Lions lemma (see Lemma I.1 p.231, in \cite{LIONS}),  $u_i\to 0$ in $L^p(\mathbb{R}^N)$, which  is impossible since $\|u_i\|_{L^p}=1$, $i=1,\ldots $. 
Thus, passing to a subsequence if necessary, we may assume that there exists $(y_i) \subset \mathbb{R}^N$ such that there holds
$
\int_{B(y_i;r)} |u_i|^q \,dx>\delta/2$, $\forall i=1,\ldots$. 
Hence, redefining  $u_i:=u_i(\cdot+y_n)$ if necessary, we have
\begin{equation*}
	\int_{B(0;r)} |u_i|^q \,dx>\delta/2, ~i=1,\ldots,
\end{equation*}
which implies that $\hat{u}^*\neq 0$.

 Recall the Brezis-Lieb lemma (see \cite{LIONS})
\begin{lem}\label{BL}
Assume $(u_n)$ is bounded in $L^\gamma(\mathbb{R}^N)$, $1\leq \gamma<+\infty$
	and $u_n \to u$ a. e. on $\mathbb{R}^N$, then
	\begin{equation*}
	\lim_{n\to +\infty}\|u_n\|_{L^\gamma}^\gamma= \|u\|_{L^\gamma}^\gamma+\lim_{n\to +\infty}\|u_n-u\|_{L^\gamma}^\gamma.
	\end{equation*}
\end{lem}
Hence,  we have 
\begin{align}
	&\|\nabla \hat{u}^*\|_{L^2}^2=\lim_{i\to \infty}\|\nabla {u_i}\|_{L^2}^2-\lim_{i\to \infty}\|\nabla( {u_i}-\hat{u}^*)\|_{L^2}^2,\\
	&\|\hat{u}^*\|_{L^p}^p=\lim_{i\to \infty}\|u_i\|_{L^p}^p-\lim_{i\to \infty}\|u_i-\hat{u}^*\|_{L^p}^p,\\
	&\|\hat{u}^*\|_{L^q}^q=\lim_{i\to \infty}\|u_i\|_{L^q}^q-\lim_{i\to \infty}\|u_i-\hat{u}^*\|_{L^q}^q.
\end{align}
Suppose that $\lim_{i\to \infty}\|\nabla( {u_i}-\hat{u}^*)\|_{L^2}^2> 0$. Then,  without loss of generality, we may assume that there holds  $\lim_{i\to \infty} \|u_i-\hat{u}^*\|_{L^q}^q>0$, $\lim_{i\to \infty}\|u_i-\hat{u}^*\|_{L^p}^p>0$ as well. Hence,
\begin{align}\label{limZer}
	\bar{\mu}=&\lim_{i\to \infty} \|\nabla {u_i}\|_{L^2}^2=\|\nabla \hat{u}^*\|_{L^2}^2+\lim_{i\to \infty}\|\nabla( {u_i}-\hat{u}^*)\|_{L^2}^2\geq \nonumber\\
	&\bar{\mu}\left(\frac{\|\hat{u}^*\|_{L^p}^\rho}{\|\hat{u}^*\|_{L^q}^\beta} +\lim_{i\to \infty}\frac{\|u_i-\hat{u}^*\|_{L^p}^\rho}{\|u_i-\hat{u}^*\|_{L^q}^\beta}\right)=\bar{\mu}\left(\frac{\|\hat{u}^*\|_{L^p}^\rho}{\|\hat{u}^*\|_{L^q}^\beta} +\frac{(1-\|\hat{u}^*\|_{L^p}^\rho)}{(1-\|\hat{u}^*\|_{L^q}^\beta)}\right)> \bar{\mu}, 
\end{align}
which  implies a contradiction. Hence $\|\nabla \hat{u}^*\|_{L^2}^2=\lim_{i\to \infty}\|\nabla u_i\|_{L^2}^2=\bar{\mu}$, $\|\hat{u}^*\|_{L^p}^p=\lim_{i\to \infty}\|u_i\|_{L^p}^p$, $\|\hat{u}^*\|_{L^q}^q=\lim_{i\to \infty}\|u_i\|_{L^q}^q$, and thus $\hat{u}^*$ is a minimizer of \eqref{MUP}.

Due to the homogeneity of $\mu^E(u)$,  we can find a minimizer  $\hat{u}^E_{\hat{\mu}^E} \in \mathcal{D}$ of \eqref{MUP} which satisfies $\sigma^E(\hat{u}^E_{\hat{\mu}^E})=1$ and $s^E(\hat{u}^E_{\hat{\mu}^E})=1$. From this it  follows that $\hat{u}^E_{\hat{\mu}^E}$ is weak solution of  \eqref{1S} with  $\mu=\hat{\mu}^E$ and prescribed energy $E_{\hat{\mu}^E}(u)=E$. Since $\mu^E(|u|)=\mu^E(u)$ for $u\in \mathcal{D} \setminus 0$, we may assume that $ \hat{u}_\mu^S\geq 0$ in $\mathbb{R}^N$. The Brezis \& Kato Theorem \cite{BK} and the $L^\gamma$ estimates for the elliptic problems  \cite{ADN} yield that $\hat{u}_\mu^E \in W^{2,\gamma}_{loc}(\mathbb{R}^N)$, for any $\gamma \in (1,+\infty)$, and whence by the regularity theory of the solutions of the elliptic problems, $\hat{u}_\mu \in C^{2}(\mathbb{R}^N)$. Consequently, the Harnack inequality \cite{trudin} implies that $\hat{u}_\mu^E> 0$ in $\mathbb{R}^N$.

\end{proof}
\begin{rem}
	The existence of spherically symmetric  ground states of the "zero-mass" case problem \eqref{1S} including  more general form
$$
-\Delta u= g(u), ~~u \in \mathcal{D}^{1,2}
$$	
was proved by Berestycki \& Lions in \cite{BerestL}. This result is obtained in \cite{BerestL} under certain hypothesis including the following 
sufficient condition
\begin{equation}\label{nesse}
	\limsup_{s\to 0+}\frac{g(s)}{s^{2^*-1}}\leq 0.
\end{equation}
Notice that in the case $2<p<q<2^*$, $\lim_{s\to 0+}g(s)/s^{2^*-1}=+\infty$.
Thus, Corollary \ref{NessC} yields that condition \eqref{nesse} is also necessary.
\end{rem}

\medskip
\bibliographystyle{amsplain}

\end{document}